\documentclass[10pt,twoside]{siamart1116}
\usepackage{amssymb}
\usepackage{eucal}
\usepackage{amsmath}
\usepackage{amscd}
\usepackage[english]{babel}
\usepackage{graphicx,epstopdf,epsfig}
\usepackage{amsfonts,epsfig,fancyhdr,graphics, hyperref,amsmath,amssymb}

\newcommand{\Title}{Classification of pre-Jordan algebras and Rota-Baxter operators on Jordan algebras in low dimensions}

\newcommand{\reals}{\mathbb{R}}
\newcommand{\complex}{\mathbb{C}}

\renewcommand{\theequation}{\thesection.\arabic{equation}}
\renewtheorem{theorem}{Theorem}[section]
\newtheorem{remark}[theorem]{Remark}
\newtheorem{example}[theorem]{Example}
\renewtheorem{lemma}[theorem]{Lemma}
\renewtheorem{corollary}[theorem]{Corollary}
\renewtheorem{proposition}[theorem]{Proposition}
\theoremstyle{definition}
\renewtheorem{definition}[theorem]{Definition}

\begin{document}
\newcommand {\emptycomment}[1]{} %to remove paragraphs

\newcommand{\nc}{\newcommand}
\newcommand{\delete}[1]{}

%%%%%%%%%%%%%%%%%%%%%%%% Statements

\renewcommand{\labelenumi}{{\rm(\alph{enumi})}}
\renewcommand{\theenumi}{\alph{enumi}}
\renewcommand{\labelenumii}{{\rm(\roman{enumii})}}
\renewcommand{\theenumii}{\roman{enumii}}

\bibliographystyle{plain}

\setcounter{page}{1}

\thispagestyle{empty}

%%%%%%%%%%%%%%%%%%%%%%%%%%%%%%%%%%%%%%%%%%%%%%%%%%%%%%%%%%%%%%%%%%

\title{\Title}

\author{Yuze Sun\thanks{School of Mathematical Sciences, Nankai University,
		Tianjin 300071, P.R.China (1810080@mail.nankai.edu.cn, \newline 1811123@mail.nankai.edu.cn,\newline 1810135@mail.nankai.edu.cn, \newline 1810278@mail.nankai.edu.cn).}
	\and Zhen Huang\footnotemark[1]
    \and Shilong Zhao\footnotemark[1]
	\and Zheshuai Tian\footnotemark[1]
}

\maketitle

\begin{abstract}
This paper is devoted to the classification of complex pre-Jordan algebras in the sense of isomorphisms in dimensions $\leq$ 3. All Rota-Baxter operators on complex Jordan algebras in dimensions $\leq$ 3 and the induced pre-Jordan algebras are also presented.

\end{abstract}



\section{Introduction}
Jordan algebras are a class of nonassociative algebras introduced in the context of
axiomatic quantum mechanics (\cite{JNW}).
Such structures were widely studied and then applied in a lot of
fields in mathematics  such as differential geometry (\cite{Chu-C}, \cite{I}, \cite{KZ}, \cite{Ko2}, \cite{Loo}), Lie theory
(\cite{J1}, \cite{Ko1}), analysis (\cite{KZ}, \cite{U})
and some fields in physics like quantum mechanics and field theory (\cite{LRH}).

Pre-Jordan algebras were  introduced in \cite{HNB} as a class of
nonassociative algebras which are closely related to Jordan
algebras. There are several motivations to study such structures.
For example, it is the algebraic structure
behind the Jordan Yang-Baxter equation which is an analogue of the
classical Yang-Baxter equation in a Lie algebra (\cite{Z2}, \cite{Z4}).
And the solutions of Jordan Yang-Baxter equation have connection with the symplectic structures on
pseudo-euclidean Jordan algebras (\cite{BB}). Given a pseudo-euclidean unital Jordan algebra, we can apply
the Tits Kantor Koecher construction (TKK construction) to obtain a quadratic Lie algebra (\cite{Ko1}, \cite{Ko3}, \cite{T}),
which is particularly important in conformal field theory because it admits a Sugawara construction (\cite{FS}).
They can be regarded as the ``opposite" of pre-Lie algebras
(\cite{Bu})  which is analogous to the Jordan algebras as the
``opposite" of Lie algebras and therefore there is a close
relationship with dendriform algebras (\cite{Lod1}) which is the
``opposite" of the relationship between pre-Lie algebras and
dendriform algebras (see the commutative diagram (4.15) in
\cite{HNB}). The operad of pre-Jordan algebras is the bisuccessor
of the operad of Jordan algebras (\cite{BBGN}), which coincides with
the fact that a Rota-Baxter operator (of weight $0$) on a Jordan
algebra induces a pre-Jordan algebra (\cite{HNB}).

Unfortunately, we have known little on pre-Jordan algebras since
we have had few explicit examples. Even except for associative
algebras which automatically satisfy the pre-Jordan axioms, we have not
known an explicit example of pre-Jordan algebras which are not associative. Simultaneously, the classification in the sense of
isomorphism is always one of the most important problems in
studying an algebraic system. Hence it is necessary and important
to classify pre-Jordan algebras in low dimensions, which is the
first motivation and aim of this paper. A direct consequence is
that not all of 2-dimensional pre-Jordan algebras are associative.

The equations involving structural constants of both Jordan and
pre-Jordan algebras are ``cubic". In general, it is difficult to
give all solutions of these equations as well as give the
corresponding classification in the sense of isomorphism, even
more for pre-Jordan algebras since they involve two identities. An
important observation is that the two identities defining
pre-Jordan algebras hold if and only if one of them holds and the
anticommutators are Jordan algebras (see
Corollary~\ref{cor:J-PJ}). So for a fixed Jordan algebra, it is
enough to solve the equations involving this identity in order to
give all compatible pre-Jordan algebras on this Jordan algebra.
Therefore there is a ``strategy" to classify pre-Jordan algebras
based on the classification of Jordan algebras of the same dimensions, which is available
in low dimensions.

Besides,  there is a closely related topic, namely
Rota-Baxter operators on Jordan algebras. Rota-Baxter operators
(on associative algebras) were introduced by G. Baxter
(\cite{Bax}) in 1960 and then they play important roles in various
areas of mathematics and mathematical physics (\cite{Guo}, \cite{R1}, \cite{R2}). The
Rota-Baxter operators on Jordan algebras also play the similar
roles and hence it is also necessary and important to give
Rota-Baxter operators on Jordan algebras explicitly.
Moreover, since a Rota-Baxter operator (of weight $0$) on a
Jordan algebra induces a (not necessarily compatible) pre-Jordan algebra, it is natural to ask
whether all pre-Jordan algebras are induced by Rota-Baxter
operators on Jordan algebras and if the answer is no, which kinds of
pre-Jordan algebras can be obtained this way.

In this paper, we commence to classify complex pre-Jordan algebras
and give Rota-Baxter operators (of weight zero) on complex Jordan
algebras in dimensions $\leq 3$. For the former, we use the
aforementioned strategy based on the known classification of
complex Jordan algebras up to dimension 3 (\cite{KS}, \cite{Wang}) and for
the latter, we find all Rota-Baxter operators on these classified
Jordan algebras and the induced pre-Jordan algebras. Comparing
these two sets of pre-Jordan algebras, we find that the  pre-Jordan
algebras from the latter only ``occupy" a small part of the former (see Theorem~\ref{2dimcompare} and Theorem~\ref{3dimcompare})
and thus we answer the above problem. Both these results on pre-Jordan algebras and Rota-Baxter operators on Jordan
algebras can be regarded as a guide for a further development.

The paper is organized as follows. In Section 2, we recall some
necessary notions and basic results on pre-Jordan algebras and
then introduce the strategy to classify compatible pre-Jordan
algebras on Jordan algebras. In Section 3, the classification of
complex pre-Jordan algebras in dimensions 1 and 2 is given through
the strategy in the previous section. In Section 4, the
classification of complex pre-Jordan algebras in dimension 3 is
given. In Section 5, we give all Rota-Baxter operators (of weight $0$) on the
Jordan algebras in dimensions $\leq 3$ and the induced pre-Jordan
algebras.

Throughout the paper, all vector spaces and algebras are over the
complex field $\complex$. We use the following notations.

Let $A$ be an $N$-dimensional vector space with a bilinear multiplication denoted
by $(x,y)\mapsto x*y$. Let $\{e_1,\cdots,e_N\}$ be a basis of
$A$. Set
\begin{equation}
e_i*e_j=\sum_{i=1}^N a_{ij}^ke_k,\;\;i,j=1,\cdots, N.
\end{equation}
Define the {\bf formal characteristic matrix} associated to
$(A,*)$ under the basis $\{e_i\}$ as
\begin{equation}
\frak M(A)=\left(\begin{matrix} e_1*e_1 & \cdots & e_1*e_N \cr
 \cdots &\cdots &\cdots\cr e_N*e_1&\cdots &
 e_N*e_N\cr\end{matrix}\right)=\left(\begin{matrix} \sum_{i=1}^N a_{11}^ke_k & \cdots & \sum_{i=1}^N a_{1N}^ke_k \cr
 \cdots &\cdots &\cdots\cr \sum_{i=1}^N a_{N1}^ke_k&\cdots &
 \sum_{i=1}^N a_{NN}^ke_k\cr\end{matrix}\right).
\end{equation}
Obviously $(A,*)$ is exactly presented by $\frak M(A)$. In this
paper, we use $\frak M(A)$ to denote $(A,*)$. The constants $a_{ij}^k$ are called the {\bf structural constants} of $(A, *)$ under the basis $\{e_i\}$.

\section{A strategy on the classification of pre-Jordan algebras}

\begin{definition}
{\rm A vector space $J$ with a bilinear multiplication $\circ : J\times
J\rightarrow J$ denoted by $(x,y)\mapsto x\circ y$ is called a
{\em Jordan algebra} if the following identities hold:
\begin{equation}x\circ y=y\circ x,\;\;\end{equation}
\begin{equation}\label{eq:JI}
((x\circ x)\circ y)\circ x = (x\circ x)\circ(y\circ
x),\end{equation} for all $x,y\in J$. A {\em homomorphism} from a
Jordan algebra $(J_1,\circ_1)$ to another Jordan algebra
$(J_2,\circ_2)$ is a linear map $\phi:J_1\rightarrow J_2$
satisfying
\begin{equation}
\phi(x\circ_1 y)=\phi(x)\circ_2 \phi(y),\;\;\forall x,y\in J_1.
\end{equation}
A bijective homomorphism between two Jordan algebras is called an
{\em isomorphism}. An isomorphism from a Jordan algebra
$(J,\circ)$ to itself is called an {\em automorphism}. The group of all
automorphisms of a Jordan algebra $(J,\circ)$ is denoted by ${\rm
Aut}(J,\circ)$.}
\end{definition}

\begin{remark}\label{rmk:Jordan}
When the characteristic of the base field is neither 2 nor 3, it
was pointed out in \cite{Al1} that for a commutative bilinear
multiplication $\circ : J\times J\rightarrow J$, Eq.~(\ref{eq:JI})
is equivalent to the following identity
\begin{equation}
(x\circ y, u, z)_\circ+(y\circ z, u, x)_\circ+(z\circ x, u,
y)_\circ=0,\;\;\forall x, y, z, u\in J,\end{equation} where $(x,
y, z)_\circ=(x\circ y)\circ z-x\circ (y\circ z)$ is the associator
of $(J,\circ)$.
\end{remark}

\begin{definition}\label{preJ} {\rm(\cite{HNB})
A vector space $A$ with a bilinear multiplication $\cdot : A\times
A\rightarrow A$ denoted by $(x,y)\mapsto x\cdot y$ is called a
{\em pre-Jordan algebra} if the following identities hold:
\begin{eqnarray}
          &&(x\circ y)\cdot(z\cdot u) + (y\circ z)\cdot(x\cdot u) + (z\circ x)\cdot(y\cdot u) \nonumber\\
          &&= x\cdot(y\cdot(z\cdot u)) + z\cdot(y\cdot(x\cdot u)) + ((x\circ z)\circ y)\cdot u \nonumber\\
          &&= z\cdot((x\circ y)\cdot u) + x\cdot((y\circ z)\cdot u) + y\cdot((z\circ x)\cdot
          u),\forall x,y,z,u\in A,
\end{eqnarray}
where $x\circ y := x\cdot y+y\cdot x$.  A
{\em homomorphism} from a pre-Jordan algebra $(A_1,\cdot_1)$ to
another pre-Jordan algebra $(A_2,\cdot_2)$ is a linear map
$\phi:A_1\rightarrow A_2$ satisfying
\begin{equation}
\phi(x\cdot_1 y)=\phi(x)\cdot_2 \phi(y),\;\;\forall x,y\in A_1.
\end{equation}
A bijective homomorphism between two pre-Jordan algebras is called
an {\em isomorphism}.}
\end{definition}

\begin{example}
According to \cite{HNB}, any associative algebra is a pre-Jordan
algebra. And obviously, the direct sum of two
pre-Jordan algebras is a pre-Jordan algebra.
\end{example}

Let $A$ be a vector space with a bilinear multiplication $\cdot :
A\times A\rightarrow A$. Let  $\circ: A\times A\rightarrow A$ be the anticommutator defined by Eq.~(\ref{eq:circ}). For any
$x,y,z,u\in J$, set
    \begin{eqnarray*}
        F(x,y,z,u) &=& (x\circ y)\cdot(z\cdot u) + (y\circ z)\cdot(x\cdot u) + (z\circ x)\cdot(y\cdot u), \\
        G(x,y,z,u) &=& x\cdot(y\cdot(z\cdot u)) + z\cdot(y\cdot(x\cdot u)) + ((x\circ z)\circ y)\cdot u, \\
        H(x,y,z,u) &=& z\cdot((x\circ y)\cdot u) + x\cdot((y\circ z)\cdot u) + y\cdot((z\circ x)\cdot u),\\
        P(x,y,z,u) &=& F(x,y,z,u) - H(x,y,z,u),\\
        Q(x,y,z,u) &=& F(x,y,z,u) - G(x,y,z,u).
\end{eqnarray*}

\begin{lemma}\label{lem:pJ}
The notations are as above.
\begin{enumerate}
\item The following identities hold:
\begin{equation}Q(x,y,z,u)=Q(z,y,x,u),\label{eq:eq1}\end{equation}
\begin{eqnarray}
&&(x\circ y,u,z)_\circ+(y\circ z,u,x)_\circ+ (z\circ x,u,y)_\circ
\nonumber\\
&&=-P(x,y,z,u) - Q(x,u,y,z) - Q(y,u,z,x) -
Q(z,u,x,y),\label{eq2}
\end{eqnarray}
for all  $x,y,z,u\in A$. \item $(A,\cdot)$ is a pre-Jordan algebra
if and only if
\begin{equation}\label{eq:PQ0}
P(x,y,z,u)=Q(x,y,z,u)=0,\;\;\forall x,y,z,u\in A. \end{equation}
\end{enumerate}
\end{lemma}

\begin{proof}
(a). Eq.~(\ref{eq:eq1}) follows from the definition of
$Q(x,y,z,u)$ immediately. Let $x,y,z,u\in A$. Then we have
\begin{eqnarray*}
&&(x\circ y,u,z)_\circ+(y\circ z,u,x)_\circ+ (z\circ
x,u,y)_\circ\\
=&& ((x\circ y)\circ u)\circ z - (x\circ y)\circ(u\circ z) + ((y\circ z)\circ u)\circ x \\
&&- (y\circ z)\circ(u\circ x) + ((z\circ x)\circ u)\circ y - (z\circ x)\circ(u\circ y)\\
=&& ((x\circ y)\circ u)\cdot z + z\cdot((x\circ y)\cdot u) + z\cdot(u\cdot(x\cdot y)) + z\cdot(u\cdot(y\cdot x))\\
&&- (x\circ y)\cdot(u\cdot z) - (x\circ y)\cdot(z\cdot u) - (u\circ z)\cdot(x\cdot y) - (u\circ z)\cdot(y\cdot x)\\
&&+ ((y\circ z)\circ u)\cdot x + x\cdot((y\circ z)\cdot u) + x\cdot(u\cdot(y\cdot z)) + x\cdot(u\cdot(z\cdot y))\\
&&- (y\circ z)\cdot(u\cdot x) - (y\circ z)\cdot(x\cdot u) - (u\circ x)\cdot(y\cdot z) - (u\circ x)\cdot(z\cdot y)\\
&&+ ((z\circ x)\circ u)\cdot y + y\cdot((z\circ x)\cdot u) + y\cdot(u\cdot(z\cdot x)) + y\cdot(u\cdot(x\cdot z))\\
&&- (z\circ x)\cdot(u\cdot y) - (z\circ x)\cdot(y\cdot u) - (u\circ y)\cdot(z\cdot x) - (u\circ y)\cdot(x\cdot z)\\
=&& -F(x,y,z,u) - F(x,u,y,z) - F(z,u,x,y) - F(y,u,z,x) \\
&&+ H(x,y,z,u) + G(x,u,y,z) + G(z,u,x,y) + G(y,u,z,x)\\
=&& -P(x,y,z,u) - Q(x,u,y,z) - Q(y,u,z,x) -Q(z,u,x,y).
\end{eqnarray*}
Hence Eq.~(\ref{eq2}) holds.

(b). It follows from Definition~\ref{preJ}.
\end{proof}

\begin{corollary}\label{cor:J-PJ}
Let $A$ be a vector space with a bilinear multiplication $\cdot :
A\times A\rightarrow A$. Let  $\circ: A\times A\rightarrow A$ be the anticommutator defined by Eq.~(\ref{eq:circ}). If
$Q(x,y,z,u)=0$ for all $x,y,z,u\in A$, then $(A,\cdot)$ is a
pre-Jordan algebra if and only if $(A,\circ)$ is a Jordan algebra.
\end{corollary}

\begin{proof}
If $(A,\cdot)$ is a pre-Jordan algebra, then $(A,\circ)$ is a
Jordan algebra by Lemma~\ref{lem:pJ} and Remark~\ref{rmk:Jordan}. Conversely, if
$(A,\circ)$ is a Jordan algebra, then by Remark~\ref{rmk:Jordan}
and Eq.~(\ref{eq2}), we have
$$P(x,y,z,u) + Q(y,z,u,x) + Q(z,u,x,y) +
Q(u,x,y,z)=0,\;\;\forall x,y,z,u\in A.$$ By assumption, we show
that $P(x,y,z,u)=0$ for all $x,y,z,u\in A$. Therefore by
Lemma~\ref{lem:pJ} (b), $(A,\cdot)$ is a pre-Jordan algebra.
\end{proof}

\begin{corollary}\label{cor:pJ-J}{\rm (\cite{HNB})} Let $(A, \cdot)$ be a  pre-Jordan algebra. Then the
anticommutator given by
\begin{equation}x\circ y=x\cdot y+y\cdot x,\quad\forall x, y\in A,\label{eq:circ} \end{equation}
defines a Jordan algebra $(J(A),\circ)$, which is called the associated
Jordan algebra of $(A,\cdot)$ and $(A,\cdot)$ is also called a
compatible pre-Jordan algebra structure on the Jordan algebra
$(J(A),\circ)$.
\end{corollary}

With the notations as above, let $\{e_1,\cdots,e_N\}$ be a basis
of $A$ and $\{e^1,\cdots,e^N\}$ be the dual basis. Set
\begin{equation}
Q^i_{klmn} = \langle
e^i,Q(e_k,e_l,e_m,e_n)\rangle,i,k,l,m,n=1,\cdots,N,
\end{equation}
where $\langle-,-\rangle$ is the usual pairing between $A$ and its dual
space $A^*$. Obviously, Eq.~(\ref{eq:PQ0}) holds if and only if
$P^i_{klmn}=Q^i_{klmn}=0$ for all $i,k,l,m,n$.

\begin{corollary}\label{variety}
Let $(J,\circ)$ be a Jordan algebra of dimension $N$. The
compatible pre-Jordan algebra structures on $J$ are 1-1
corresponding to the common roots of a set of
polynomials. More precisely, the set contains at most
$\frac{N^4(N+1)}{2}$ polynomials of degree at most $3$ with
$\frac{N^2(N-1)}{2}$ indeterminates.
\end{corollary}
\begin{proof} Let $\{e_1,\cdots,e_N\}$ be a  basis of $J$. Set
$$e_i\circ e_j=\sum_{k=1}^{N}c_{ij}^ke_k, \;\;i,j=1,\cdots, N.$$
That is, $c_{ij}^k$ are the structural constants of the Jordan
algebra $(J,\circ)$. Let $(J,\cdot)$ be a compatible pre-Jordan
algebra on $(J,\circ)$. Then by Eq.~(\ref{eq:circ}), we can assume
that \begin{equation}\label{eq:assume} e_i\cdot e_j=\begin{cases}
        \sum_{k=1}^{N}(\frac{c_{ij}^k}{2}+x_{ij}^k)e_k, & \mbox{if $i<j$} \\
        \sum_{k=1}^{N}\frac{c_{ij}^k}{2}e_k, & \mbox{if $i=j$} \\
        \sum_{k=1}^{N}(\frac{c_{ij}^k}{2}-x_{ji}^k)e_k, & \mbox{if
        $i>j$},
      \end{cases}\end{equation}
where $x_{ij}^k$ are indeterminates, $1\leq i<j\leq n,1\leq k\leq
n$. Note that the number of these indeterminates is
$\frac{N^2(N-1)}{2}$. Furthermore, by Corollary~\ref{cor:J-PJ},
$(J,\cdot)$ is a pre-Jordan algebra if and only if
$$Q^i_{klmn}=0,\forall i,k,l,m,n=1,\cdots, N.$$\\
If we set $x_{ij}^k=-x_{ji}^k$ for $i> j$ and $x_{ii}^k=0$, we can write the equations as:\\
\begin{eqnarray*}
  Q^i_{klmn} &=& \sum_{r=1}^{N}\sum_{j=1}^{N}\left\{\left[\left(\frac{c^j_{mn}}{2}+x_{mn}^j\right)c_{kl}^r+\left(\frac{c_{kn}^j}{2}+x_{kn}^j\right)c_{lm}^r+\left(\frac{c_{ln}^j}{2}+x_{ln}^j\right)c_{mk}^r\right]\left(\frac{c^i_{rj}}{2}+x^i_{rj}\right) \right.\\
   &+& \left(\frac{c^r_{mn}}{2}+x_{mn}^r\right)\left(\frac{c_{lr}^j}{2}+x_{lr}^j\right)\left(\frac{c_{kj}^i}{2}+x_{kj}^i\right)+\left(\frac{c_{kn}^r}{2}+x_{kn}^r\right)\left(\frac{c_{lr}^j}{2}+x_{lr}^j\right)\left(\frac{c_{mj}^i}{2}+x_{mj}^i\right)\\
&+&\left.c_{km}^rc_{rl}^j\left(\frac{c_{jn}^i}{2}+x_{jn}^i\right)\right\}\\
   &=&0.
\end{eqnarray*}
By Eq.~(\ref{eq:eq1}), we can assume that $k\leq m$. Thus the
number of these equations with indeterminates $x_{ij}^k$ is at
most $\frac{N^4(N+1)}{2}$. Therefore the conclusion follows.
\end{proof}

\begin{lemma}\label{iso} Let $(J_1,\circ_1)$ and $(J_2,\circ_2)$ be
two Jordan algebras and $\phi:(J_1,\circ_1)\rightarrow (J_2,\circ_2)$ be
an isomorphism. If $(J_1,\cdot_1)$ is a compatible pre-Jordan
algebra structure on $(J_1,\circ_1)$, then the bilinear
multiplication $\cdot_2:J_2\times J_2\rightarrow J_2$ defined by
\begin{equation}\label{eq:iso}
x\cdot_2y:=\phi(\phi^{-1}(x)\cdot_1\phi^{-1}(y)),\;\;\forall
x,y\in J_2,
\end{equation}
gives a compatible pre-Jordan algebra structure on $(J_2,\circ_2)$
which is isomorphic to $(J_1,\cdot_1)$.
\end{lemma}

\begin{proof}
It is straightforward to show that $(J_2,\cdot_2)$ is a pre-Jordan
algebra which is isomorphic to $(J_1,\cdot_1)$. Moreover, for any
$x,y\in J_2$, we have
$$x\cdot_2y+y\cdot_2x=\phi(\phi^{-1}(x)\cdot_1\phi^{-1}(y))+\phi(\phi^{-1}(y)\cdot_1\phi^{-1}(x))=\phi(\phi^{-1}(x)\circ_1\phi^{-1}(y))=x\circ_2y.$$
Therefore $(J_2,\cdot_2)$ is a  compatible pre-Jordan algebra
structure on $(J_2,\circ_2)$.
\end{proof}

\begin{corollary}\label{orbit}
Let $(J,\circ)$ be a Jordan algebra. Let $(J,\cdot_1)$ and
$(J,\cdot_2)$ be two compatible pre-Jordan algebra structures on
$(J,\circ)$. Then $(J,\cdot_1)$ and $(J,\cdot_2)$ are isomorphic
and $\phi:(J,\cdot_1)\rightarrow (J,\cdot_2)$ is an isomorphism if and only if
$\phi\in {\rm Aut}(J,\circ)$ and $\phi$ satisfies
Eq.~(\ref{eq:iso}). Equivalently,  ${\rm Aut}(J,\circ)$ acts on
the set of all compatible pre-Jordan algebra structures on
$(J,\circ)$ through Eq.~(\ref{eq:iso}), and the isomorphism
classes of these pre-Jordan algebras are exactly the orbits of
this action.
\end{corollary}

\begin{proof} ``$\Rightarrow$": If $\phi:J\rightarrow J$
is an isomorphism of pre-Jordan algebras from $(J,\cdot_1)$ to
$(J,\cdot_2)$, then it is straightforward to show that $\phi$ is
an isomorphism from $(J,\circ)$ to itself, that is, $\phi\in {\rm Aut}(J,\circ)$. Moreover,
$\phi$ satisfies Eq.~(\ref{eq:iso}) by Definition~\ref{preJ}.

``$\Leftarrow$": It follows from Lemma~\ref{iso} by letting
$(J_1,\circ_1)=(J_2,\circ_2)=(J,\circ)$.
\end{proof}

Therefore, we give our strategy on classifying finite-dimensional
complex pre-Jordan algebras in the sense of isomorphism as
follows. We divide it into several steps.

{\bf Step 1} We classify finite-dimensional complex Jordan
algebras in the sense of isomorphism.

{\bf Step 2} For a fixed Jordan algebra $(J,\circ)$ in Step
1, that is, with a fixed basis and fixed structural constants
$c_{ij}^k$, we give all solutions of the set of equations
$\{Q^i_{klmn}=0\}$ with the indeterminates $x_{ij}^k$. Thus we get
all compatible pre-Jordan algebras on $(J,\circ)$.

{\bf Step 3} For the fixed Jordan algebra $(J,\circ)$ in
Step 2, we first give ${\rm Aut}(J,\circ)$ and then determine the
orbits by its action on the set of compatible pre-Jordan algebra
structures (in fact, the solutions of the set of equations
$\{Q^i_{klmn}=0\}$) in Step 2 through Eq.~(\ref{eq:iso}).

In the following sections, we will illustrate that such a strategy
is available, in particular in low dimensions. Note that  the
classification of complex Jordan algebras in dimensions $\leq 3$
is known, that is, Step 1 has been already finished in
our cases.

\section{Classification of complex pre-Jordan algebras in dimensions $\leq2$}

In this section, we  classify all complex pre-Jordan algebras in
the sense of isomorphisms in dimensions up to 2.

\begin{proposition}\label{prop:1dim}
Let $\{e\}$ be a basis of a 1-dimensional vector space $A$. Then
there are exactly 2 non-isomorphic pre-Jordan algebras
$(A,\cdot_1)$ and $(A,\cdot_2)$ on $A$ given by $e\cdot e=e$ and
$e\cdot e=0$ respectively. Both of them are associative.
\end{proposition}

\begin{proof}
In fact, we have $e\cdot e=\lambda e$, where $\lambda\in \complex$. If $\lambda\ne 0$, then by a linear transformation given by
$e\mapsto \lambda e'$, we get $e'\cdot e'=e'$. It is obvious
that both of them are associative and hence they are pre-Jordan
algebras.
\end{proof}

\begin{lemma}{\rm (\cite{KS})}
    Every 2-dimensional Jordan algebra is isomorphic to one of the following (mutually non-isomorphic) Jordan
    algebras given by their formal characteristic matrices respectively.

    $$J_1:\begin{pmatrix}
     2e_1  &  0   \\
      0 & 2e_2
\end{pmatrix}, J_2:\begin{pmatrix}
2e_1  &  2e_2   \\
2e_2 & 0
\end{pmatrix},
 J_3 :\begin{pmatrix}
2e_1  &  0   \\
0 & 0
\end{pmatrix},$$

 $$J_4:\begin{pmatrix}
2e_1  &  e_2   \\
e_2 & 0
\end{pmatrix},
J_5:\begin{pmatrix}
2e_2  &  0   \\
0 & 0
\end{pmatrix},
     J_6:\begin{pmatrix}
0  &  0   \\
0 & 0
\end{pmatrix}.$$
\end{lemma}

In the following content, we use the notations given in the proof of Corollary~\ref{variety}
(see Eq.~(\ref{eq:assume})).
As an illustration, we give an explicit proof for the
classification of the compatible pre-Jordan algebras on the Jordan
algebra $(J_1,\circ)$, whereas the proofs for the other cases are
shortened by giving the ``essentially available " equations only.

\begin{proposition}
    There are exactly two non-isomorphic compatible pre-Jordan algebra structures on  the Jordan algebra $(J_1,\circ)$ given by the following formal characteristic matrices respectively.
$$J_{1,1}:\begin{pmatrix}
e_1  &  0   \\
0 & e_2
\end{pmatrix},\quad
     J_{1, 2}:\begin{pmatrix}
e_1  &  e_2   \\
-e_2 & e_2

\end{pmatrix}.$$

\end{proposition}

\begin{proof}
The set of equations $\{Q^i_{jklm}=0\}$ is written as follows.
\begin{eqnarray*}
&&Q^{1}_{1111}=Q_{1111}^2=0,\;\;
Q_{1112}^1=-2x_{12}^1x_{12}^2(x_{12}^2-2)=0,\;\;
\\&&
Q_{1112}^2=-2x_{12}^2(x_{12}^2-2)(x_{12}^2-1)=0,\;\;
Q_{1121}^1=x_{12}^1x_{12}^2(x_{12}^2-1)=0,\;\;
\\&&
Q_{1121}^2=x_{12}^2(x_{12}^2-1)^2=0,\;\;
Q_{1122}^1=x_{12}^1(x_{12}^1x_{12}^2+x_{12}^1-x_{12}^2+1)=0,\;\;
\\&&
Q_{1122}^2=x_{12}^2(x_{12}^1x_{12}^2+x_{12}^1-2x_{12}^2+2)=0,\;\;
Q_{1211}^1=Q_{1211}^2=0,\;\;
\\&&
Q_{1212}^1=2x_{12}^1(x_{12}^1x_{12}^2+x_{12}^1-x_{12}^2+1)=0,\;\;
Q_{1212}^2=2x_{12}^2(x_{12}^1x_{12}^2-x_{12}^2+1)=0,\;\;
\\&&
Q_{1221}^1=-x_{12}^1(x_{12}^1x_{12}^2+2x_{12}^1-x_{12}^2+2)=0,\;\;
Q_{1221}^2=-x_{12}^2(x_{12}^1x_{12}^2+x_{12}^1-x_{12}^2+1)=0,\;\;
\\&&
Q_{1222}^1=-x_{12}^1(x_{12}^1+1)^2=0,\;\;
Q_{1222}^2=-x_{12}^1x_{12}^2(x_{12}^1+1)=0,\;\;
\\&&
Q_{2121}^1=-2x_{12}^1(x_{12}^1x_{12}^2+x_{12}^1+1),\;\;
Q_{2121}^2=-2x_{12}^2(x_{12}^1x_{12}^2+x_{12}^1-x_{12}^2+1)=0,\;\;
\\&&
Q_{2122}^1=Q_{2122}^2=0,\;\;
Q_{2221}^1=2x_{12}^1(x_{12}^1+1)(x_{12}^1+2),\;\;
\\&&
Q_{2221}^2=2x_{12}^1x_{12}^2(x_{12}^1+2),\;\;
Q_{2222}^1=Q_{2222}^2=0.\;\;
\end{eqnarray*}
In fact, the
above equations hold if and only if the following (``essentially
available") equations hold:
$$Q^1_{1112}=-2x_{12}^1x_{12}^2(x_{12}^2-2)=0,\;\;
    Q^2_{1121}=x_{12}^2(x_{12}^2-1)^2=0,\;\;
    Q^1_{1222}=-x_{12}^1(x_{12}^1+1)^2=0.
$$
It is straightforward to get all solutions as follows.
$$(1)\; x_{12}^1=0, x_{12}^2=0;\;\;(2)\; x_{12}^1=0, x_{12}^2=1;\;\;(3)\;x_{12}^1=-1, x_{12}^2=0.$$
Therefore they correspond to the following pre-Jordan algebras
given by the following formal characteristic matrices
respectively.
$$J_{1,1}:\begin{pmatrix}
e_1  &  0   \\
0 & e_2
\end{pmatrix},\;\;
  J_{1,2}:\begin{pmatrix}
e_1  &  e_2   \\
-e_2 & e_2

\end{pmatrix},\;\;
   J_{1,3}:\begin{pmatrix}
e_1  &  -e_1   \\
e_1 & e_2

\end{pmatrix}.$$
Note that $(J_{1,2},\cdot)$ is isomorphic to $(J_{1,3},\cdot)$ by
a linear transformation given by $e_1\mapsto e_2,
e_2\mapsto e_1$. Moreover, $(J_{1,1},\cdot)$ is not isomorphic
to $(J_{1,2},\cdot)$ since the former is commutative and
associative, whereas the latter is neither commutative nor
associative.
\end{proof}

\begin{proposition}
     There are exactly two non-isomorphic compatible pre-Jordan algebra structures on  the Jordan algebra $(J_2,\circ)$ given by the following formal characteristic matrices respectively.
$$J_{2,1}:\begin{pmatrix}
e_1  &  e_2   \\
e_2 & 0

\end{pmatrix},\;\;J_{2,2}:\begin{pmatrix}
e_1  &  2e_2   \\
0 & 0

\end{pmatrix}
.$$
\end{proposition}

\begin{proof} The ``essentially available" equations in the set of
equations $\{Q^i_{jklm}=0\}$ are
$$Q^1_{1222}=-(x_{12}^1)^3=0,\;\;
    Q^2_{1121}=(x_{12}^2)^2(x_{12}^2-1)=0.$$
All the solutions of these equations are
$$(1)\;x_{12}^1=0,x_{12}^2=0;\;\;(2)\;x_{12}^1=0,x_{12}^2=1.$$
They correspond to the pre-Jordan algebras $(J_{2,1},\cdot)$ and
$(J_{2,2},\cdot)$ respectively. Moreover, $(J_{2,1},\cdot)$ is not
isomorphic to $(J_{2,2},\cdot)$ since the former is commutative
and associative, whereas the latter is neither commutative nor
associative.
\end{proof}

\begin{proposition}
 There are exactly two non-isomorphic compatible pre-Jordan algebra structures on the Jordan algebra $(J_3,\circ)$ given by the following formal characteristic matrices respectively.
$$J_{3,1}:\begin{pmatrix}
e_1  &  0   \\
0 & 0

\end{pmatrix},\;\;
    J_{3,2}:\begin{pmatrix}
e_1  &  e_2   \\
-e_2 & 0

\end{pmatrix}.$$

\end{proposition}

\begin{proof}
    The ``essentially available" equations in the set of
equations $\{Q^i_{jklm}=0\}$ are
 $$
    Q^1_{1222}=-(x_{12}^1)^3=0,\;\;
    Q^2_{1121}=x_{12}^2(x_{12}^2-1)^2=0.$$
All the solutions of these equations are
$$(1)\;x_{12}^1=0,x_{12}^2=0;\;\;(2)\;x_{12}^1=0,x_{12}^2=1.$$
They correspond to the pre-Jordan algebras $(J_{3,1},\cdot)$ and
$(J_{3,2},\cdot)$ respectively. Moreover, $(J_{3,1},\cdot)$ is not
isomorphic to $(J_{3,2},\cdot)$ since the former is commutative
and associative, whereas the latter is neither commutative nor
associative.
\end{proof}

\begin{proposition}
  There are exactly three non-isomorphic compatible pre-Jordan algebra structures on the Jordan algebra $(J_4,\circ)$ given by the following formal characteristic matrices respectively.
$$
    J_{4,1}:\begin{pmatrix}
e_1  &  0   \\
e_2 & 0

\end{pmatrix},\;\;
    J_{4,2}:\begin{pmatrix}
e_1  &  e_2   \\
0 & 0

\end{pmatrix},\;\;
    J_{4,3}:\begin{pmatrix}
e_1  &  2e_2   \\
-e_2 & 0

\end{pmatrix}.$$
\end{proposition}

\begin{proof}
    The ``essentially available" equations in the set of
equations $\{Q^i_{jklm}=0\}$ are
   $$ Q^1_{1222}=-(x_{12}^1)^3=0,\;\;
    8Q^2_{1121}=(2x_{12}^2-3)(2x_{12}^2-1)(2x_{12}^2+1)=0.$$
All the solutions of these equations are
$$(1)\;x_{12}^1=0,x_{12}^2=-\frac{1}{2};\;\;(2)\;x_{12}^1=0,x_{12}^2=\frac{1}{2};\;\;(3)\;x_{12}^1=0,x_{12}^2=\frac{3}{2}.$$
They correspond to the pre-Jordan algebras $(J_{4,1},\cdot)$,
$(J_{4,2},\cdot)$ and $(J_{4,3},\cdot)$ respectively. It is
straightforward to know that both $(J_{4,1},\cdot)$ and
$(J_{4,2},\cdot)$ are associative whereas $(J_{4,3},\cdot)$ is not
associative since in $(J_{4,3},\cdot)$, we have $$2e_2=(e_1\cdot
e_1)\cdot e_2\neq e_1\cdot (e_1\cdot e_2)=4e_2.$$
Moreover, $(J_{4,1},\cdot)$ is not isomorphic to $(J_{4,2},\cdot)$
since $e_1$ is a right identity in the former whereas there is not
a right identity in the latter since $e_2\cdot x=0, for all x\in J_{4,2}$.
\end{proof}

\begin{proposition}
 There is exactly one compatible pre-Jordan algebra structure on the Jordan algebra $(J_5,\circ)$ given by the following formal characteristic matrix.
$$J_{5,1}:\begin{pmatrix}
e_2  &  0   \\
0 & 0
\end{pmatrix}.$$
\end{proposition}

\begin{proof}
   The ``essentially available" equations in the set of
equations $\{Q^i_{jklm}=0\}$ are
  $$Q^1_{1211}=2(x_{12}^1)^2=0,\;\;
        Q^2_{1111}=-2(x_{12}^1+(x_{12}^2)^2)=0.$$
There is only one solution
$$a=b=0.$$
It corresponds to the pre-Jordan algebra $(J_{5,1},\cdot)$.
\end{proof}

\begin{proposition}
  There is exactly one compatible pre-Jordan algebra structure on the Jordan algebra $(J_6,\circ)$ given by the following formal characteristic matrix.
$$J_{6,1}:\begin{pmatrix}
0 &  0   \\
0 & 0
\end{pmatrix}.$$
\end{proposition}

\begin{proof}
The ``essentially available" equations in the set of equations
$\{Q^i_{jklm}=0\}$ are
  $$Q^1_{1222}=-(x_{12}^1)^3=0,\;\;
    Q^2_{1121}=(x_{12}^2)^3=0.$$
There is only one solution
$$x_{12}^1=x_{12}^2=0.$$
It corresponds to the pre-Jordan algebra $(J_{6,1},\cdot)$.
\end{proof}

Summarize the above study in this section, we have the following
conclusion.

\begin{theorem}
Every 2-dimensional pre-Jordan algebra is isomorphic to one of the
following (mutually non-isomorphic) pre-Jordan
    algebras given by their formal characteristic matrices respectively.
$$J_{1,1}:\begin{pmatrix}
e_1  &  0   \\
0 & e_2
\end{pmatrix},\;
     J_{1, 2}:\begin{pmatrix}
e_1  &  e_2   \\
-e_2 & e_2

\end{pmatrix},\;J_{2,1}:\begin{pmatrix}
e_1  &  e_2   \\
e_2 & 0

\end{pmatrix},\;J_{2,2}:\begin{pmatrix}
e_1  &  2e_2   \\
0 & 0

\end{pmatrix},\;$$
$$J_{3,1}:\begin{pmatrix}
e_1  &  0   \\
0 & 0

\end{pmatrix},\;
  J_{3,2}:\begin{pmatrix}
e_1  &  e_2   \\
-e_2 & 0

\end{pmatrix},\;
J_{4,1}:\begin{pmatrix}
e_1  &  0   \\
e_2 & 0

\end{pmatrix},\;
    J_{4,2}:\begin{pmatrix}
e_1  &  e_2   \\
0 & 0

\end{pmatrix},\;$$
  $$  J_{4,3}:\begin{pmatrix}
e_1  &  2e_2   \\
-e_2 & 0

\end{pmatrix},\;
J_{5,1}:\begin{pmatrix}
e_2  &  0   \\
0 & 0
\end{pmatrix},\;
J_{6,1}:\begin{pmatrix}
0 &  0   \\
0 & 0
\end{pmatrix}.$$
Moreover, $(J_{1,1},\cdot)$, $(J_{2,1},\cdot)$, $(J_{3,1},\cdot)$,
$(J_{5,1},\cdot)$ and $(J_{6,1},\cdot)$ are commutative and
associative, and $(J_{4,1},\cdot)$ and $(J_{4,2},\cdot)$ are
associative.
\end{theorem}

\section{Classification of 3-dimensional complex pre-Jordan algebras}

In this section, we classify 3-dimensional complex  pre-Jordan
algebras in the sense of isomorphism based on the classification
of 3-dimensional complex Jordan algebras given in \cite{KS}.

According to \cite{KS},
there are exactly 20 classes of 3-dimensional Jordan algebras:
\begin{enumerate}
\item[(1)] $\mathcal{J}_i, 1\leq i\leq 8$ are the Jordan algebras which are not associaive; \item[(2)] $A_{{\bf 1}_i}, 1\leq i\leq4$ are the
associative unitary Jordan algebras; \item[(3)] $A_{i}, 1\leq
i\leq8$ are the associative but not unitary Jordan algebras.
\end{enumerate}

As an illustration, we give
an explicit proof for the classification of the compatible
pre-Jordan algebras on the Jordan algebra $({\mathcal
J}_1,\circ)$, whereas the proofs for the other cases are omitted.

\begin{proposition} {\rm (case $\mathcal{J}_1$)} The formal
characteristic matrix of the Jordan algebra
$(\mathcal{J}_1,\circ)$ is given by $$\begin{pmatrix}
2e_1  &  2e_2  &2e_3 \\
2e_2 & 0  &2e_1\\
2e_3 & 2e_1 &0
\end{pmatrix}.$$
Any compatible pre-Jordan algebra structure on
$(\mathcal{J}_1,\circ)$ is isomorphic to one of the following
(mutually non-isomorphic) pre-Jordan
    algebras given by their formal characteristic matrices respectively.
$$\mathcal {J}_{1,1}(\alpha):
\begin{pmatrix}
e_1  &  e_1+(1+\lambda)e_2+\frac{e_3}{2\lambda}  &2\alpha e_1+2\lambda \alpha e_2+(2-\lambda)e_3 \\
-e_1+(1-\lambda)e_2-\frac{e_3}{2\lambda} & 0  &2\lambda e_1-2\alpha e_2+e_3\\
-2\alpha e_1-2\lambda \alpha e_2+\lambda e_3 & (2-2\lambda)e_1+2\alpha e_2-e_3 &0
\end{pmatrix},$$ where $\alpha \in \complex$ and $\lambda$ is a nonzero root of
$x^2-x+\alpha =0$, and different choices of $\lambda$ give isomorphic pre-Jordan algebras,\\
$$ \mathcal {J}_{1,2}:\begin{pmatrix}
e_1  &  2e_2  &e_3 \\
0  &  0  &2e_1\\
e_3  &  0  &0
\end{pmatrix},\;\;
\mathcal {J}_{1,3}:\begin{pmatrix}
e_1  &  2e_2  &e_2+e_3 \\
0  &  0  &2e_1\\
-e_2+e_3  &  0  &0
\end{pmatrix}.$$
\end{proposition}
\begin{proof}
We determine the structural constants. By Corollary
\ref{variety},
\begin{eqnarray*}
&&Q^3_{1211}+Q^3_{1121}+Q^3_{1112}=-2x_{12}^3(x_{12}^2+x_{13}^3-1)=0,\\
&&Q^2_{1311}+Q^2_{1131}+Q^2_{1113}=-2x_{13}^2(x_{12}^2+x_{13}^3-1)=0,\\
&&2Q^3_{1221}+Q^3_{1122}=2x_{12}^3(x_{12}^1-x_{23}^3)=0,\\
&&Q^2_{1321}+Q^2_{1231}+Q^2_{1123}=2x_{13}^2(x_{12}^1-x_{23}^3)=0,\\
&&Q^3_{1231}+Q^3_{1321}+Q^3_{1132}=2x_{12}^3(x_{13}^1+x_{23}^2)=0,\\
&&2Q^2_{1331}+Q^3_{1133}=2x_{13}^2(x_{13}^1+x_{23}^2)=0.
\end{eqnarray*}

If at least one of $\begin{cases}
x_{23}^3=x_{12}^1\\
x_{12}^2+x_{13}^3=1\\
x_{13}^1+x_{23}^2=0\\
\end{cases}$ fails to hold, then $x_{12}^3=x_{13}^2=0$. Thus,
\begin{equation*}
Q^1_{2221}=2(x_{12}^1)^3=0,\;\;
Q^1_{3331}=2(x_{13}^1)^3=0,\;\;
Q^3_{2232}=(x_{23}^3)^3=0,\;\;
Q^2_{2333}=-(x_{23}^2)^3=0.
\end{equation*}
And $x_{12}^2+x_{13}^3-1=\frac{1}{8}(3Q^1_{2233}+3Q^1_{2332}-2Q^3_{3132}-2Q^2_{2123})=0$.
These lead to a contradiction, so $\begin{cases}
x_{23}^3=x_{12}^1\\
x_{12}^2+x_{13}^3=1\\
x_{13}^1+x_{23}^2=0\\
\end{cases}$ always holds. What's more, $x_{23}^1-2x_{12}^2+1=\frac{1}{2}(Q^1_{1123}-Q^1_{1132})=0$.\\
We have reduced 4 indeterminates, and the remaining unknown elements follow these ``essentially available" equations:
\begin{eqnarray*}
x_{12}^1(1-x_{12}^2)=x_{12}^3x_{13}^1&\Leftrightarrow& Q^1_{1211}+Q^1_{1121}+Q^1_{1112}=0,\\
x_{12}^1x_{13}^2=x_{12}^2x_{13}^1&\Leftrightarrow& Q^1_{1311}+Q^1_{1131}+Q^1_{1113}=0,\\
x_{12}^2(1-x_{12}^2)=x_{12}^3x_{13}^2&\Leftrightarrow& Q^2_{1211}+Q^2_{1121}+Q^2_{1112}=0,\\
(x_{12}^1)^2=2x_{12}^2x_{12}^3&\Leftrightarrow& Q^1_{2322}+Q^1_{2232}+Q^1_{2223}=0,\\
(x_{13}^1)^2=2x_{13}^2(1-x_{12}^2)&\Leftrightarrow& Q^1_{3233}+Q^1_{2333}+Q^1_{3332}=0,\\
x_{12}^1x_{13}^1=2x_{12}^2x_{13}^3&\Leftrightarrow& 2Q^1_{2332}+Q^1_{2233}=0.
\end{eqnarray*}
It has been shown that the set of solutions mentioned in Corollary~\ref{variety} 1-1 corresponds to the subset of $\complex^5$ cut by the equations above. Their coordinates are $(x_{12}^1,x_{12}^2,x_{12}^3,x_{13}^1,x_{13}^2)$. By Corollary~\ref{orbit}, the next step is to investigate how ${\rm Aut}(\mathcal{J}_1, \circ)$ acts on this set. It is straightforward to show that
${\rm Aut}(\mathcal{J}_1,\circ)=H\cup\tau H$, where
$H=\{{\rm diag}(1,\mu,\mu^{-1})|\mu\neq0\}$ and
$\tau=\begin{pmatrix}
                                                                                    1 & 0 & 0 \\
                                                                                    0 & 0 & 1 \\
                                                                                    0 & 1 & 0
                                                                                  \end{pmatrix}$.\\
And by simple calculation,
\[
  \begin{aligned}
  {\rm diag}(1,\mu,\mu^{-1})(x_{12}^1,x_{12}^2,x_{12}^3,x_{13}^1,x_{13}^2)&=(x_{12}^1\mu^{-1},x_{12}^2,x_{12}^3\mu^{-2},x_{13}^1\mu,x_{13}^2\mu^2),\\
\tau(x_{12}^1,x_{12}^2,x_{12}^3,x_{13}^1,x_{13}^2)&=(x_{13}^1,1-x_{12}^2,x_{13}^2,x_{12}^1,x_{12}^3).
  \end{aligned}
\]
Notice that $\alpha =x_{12}^2(1-x_{12}^2)=x_{12}^3x_{13}^2$ is invariant under ${\rm Aut}(\mathcal{J}_1, \circ)$. And then we introduce a lemma.\\
Suppose $G$ is a group acting on a set $S$ by $\phi$ and $G^\prime$ is a normal subgroup of $G$. Note the set of all orbits by $S/G$. $G$ could act on $S/G^\prime$ by $\phi^\prime$: $\phi^\prime(g)([x])=[\phi(g)x]$ (note the orbit of $x$ in $S/G^\prime$ by $[x]$). Hence $\phi^\prime(g)=$id if and only if $g\in G^\prime$. That is, the action of $G/G^\prime$ on $S/G^\prime$ by $\pi: gG^\prime\mapsto\phi(g)$ is well defined. And $\pi(gG)[x]=[y]$ if and only if there is $h\in G^\prime$, such that $\phi(gh)(x)=y$, which tells the following lemma.
\begin{lemma}
Suppose $G$ is a group acting on a set $S$, and $G^\prime\lhd G$ is a normal subgroup. Note the orbits of $G^\prime$ acting on $S$ by $S/G^\prime$. Then $G/G^\prime$ acts on $S/G^\prime$, and $S/G\simeq(S/G^\prime)/(G/G^\prime)$.
\end{lemma}
\indent When $\alpha \neq 0$, $x_{12}^1x_{13}^1=2\alpha \neq0$. Apply the lemma above on $S=\{\mbox{tuples with } \alpha \neq0\}$, $G={\rm Aut}(\mathcal{J}_1, \circ)$ and $G^\prime=H=\mathrm{ker}[\mbox{det: }G\rightarrow\complex^*]$. For each orbit in $S/H$, we fix a representing element with $x_{12}^1=1$, and in this way we can reduce the problem to determining the orbit of the induced action of $\{1,\tau\}$ on the tuples with the first term $1$. Consider an arbitrary tuple $(1,x_{12}^2,(2x_{12}^2)^{-1},2\alpha ,2x_{12}^2\alpha )$. Apply $\tau$ on this tuple and it becomes $(1,1-x_{12}^2,(2-2x_{12}^2)^{-1},2\alpha ,2(1-x_{12}^2)\alpha )$. Noticing the fact that $x_{12}^2$ and $1-x_{12}^2$ are the roots of $x^2-x+\alpha =0$, the discussion above ensures that a different choice of the nonzero root $x_{12}^2$ of this equation does not make any difference on the isomorphism class. This proves the part $\mathcal{J}_{1,1}(\alpha \neq0)$.\\
\indent When $\alpha =0=x_{12}^2(1-x_{12}^2)$, $x_{12}^2=0$ or $1$. Because $\tau$ sends the tuples with $x_{12}^2=0$ to the ones with $x_{12}^2=1$, every orbit admits at least one element with $x_{12}^2=1$. It is obvious that $\{\mbox{tuples with }\alpha =0\}/G\simeq\{\mbox{tuples with }x_{12}^2=1\}/H$. When $x_{12}^2=1$, the equations are reduced to

\begin{equation*}
x_{13}^1=0,\;\;
(x_{12}^1)^2=2x_{12}^3,\;\;
x_{12}^1x_{13}^2=0.
\end{equation*}
If $x_{12}^1\neq0$, it can be fixed to $1$, and this is the case $\mathcal{J}_{1,1}(\alpha =0)$.\\
If $x_{12}^1=0$, then $x_{13}^2=0$ or $1$, generating $\mathcal{J}_{1,2}$ and $\mathcal{J}_{1,3}$ respectively.
\end{proof}

\begin{proposition}{\rm (case $\mathcal{J}_2$)} The formal
characteristic matrix of the Jordan algebra
$(\mathcal{J}_2,\circ)$ is given by
$$\begin{pmatrix}
2e_1  &  0  &e_3 \\
0  &  2e_2  &e_3\\
e_3  &  e_3  &0
\end{pmatrix}.$$Any compatible pre-Jordan algebra structure on
$(\mathcal{J}_2,\circ)$ is isomorphic to one of the following
(mutually non-isomorphic) pre-Jordan
algebras given by their formal characteristic matrices respectively.\begin{eqnarray*}
&&\mathcal{J}_{2,1}:\begin{pmatrix}
e_1  &  -e_1  &e_3 \\
e_1  &  e_2  &0\\
0  &  e_3  &0
\end{pmatrix},
\mathcal{J}_{2,2}:\begin{pmatrix}
e_1  &  0  &e_3 \\
0  &  e_2  &e_3\\
0  &  0  &0
\end{pmatrix},
\mathcal{J}_{2,3}:\begin{pmatrix}
e_1  &  0  &e_3 \\
0  &  e_2  &e_1+e_3\\
0  &  -e_1  &0
\end{pmatrix},\\&&
\mathcal{J}_{2,4}:\begin{pmatrix}
e_1  &  0  &0 \\
0  &  e_2  &e_3\\
e_3  &  0  &0
\end{pmatrix},
\mathcal{J}_{2,5}:\begin{pmatrix}
e_1  &  -e_1  &0 \\
e_1  &  e_2  &2e_3\\
e_3  &  -e_3  &0
\end{pmatrix},
\mathcal{J}_{2,6}:\begin{pmatrix}
e_1  &  -e_1  &e_3 \\
e_1  &  e_2  &e_3\\
0  &  0  &0
\end{pmatrix},\\&&
\mathcal{J}_{2,7}:\begin{pmatrix}
e_1  &  -e_1  &2e_3 \\
e_1  &  e_2  &0\\
-e_3  &  e_3  &0
\end{pmatrix}.
\end{eqnarray*}
\\
\end{proposition}
\begin{proposition}{\rm (case $\mathcal{J}_3$)} The formal
characteristic matrix of the Jordan algebra
$(\mathcal{J}_3,\circ)$ is given by$$\begin{pmatrix}
2e_1  &  0  &0 \\
0  &  2e_2  &e_3\\
0  &  e_3  &0
\end{pmatrix}.$$Any compatible pre-Jordan algebra structure on
$(\mathcal{J}_3,\circ)$ is isomorphic to one of the following
(mutually non-isomorphic) pre-Jordan
algebras given by their formal characteristic matrices respectively.\\
\begin{eqnarray*}
&&\mathcal{J}_{3,1}:\begin{pmatrix}
e_1  &  -e_1  &0 \\
e_1  &  e_2  &0\\
0  &  e_3  &0
\end{pmatrix},
\mathcal{J}_{3,2}:\begin{pmatrix}
e_1  &  -e_1  &0 \\
e_1  &  e_2  &e_3\\
0  &  0  &0
\end{pmatrix},
\mathcal{J}_{3,3}:\begin{pmatrix}
e_1  &  -e_1  &0 \\
e_1  &  e_2  &2e_3\\
0  &  -e_3  &0
\end{pmatrix},\\&&
\mathcal{J}_{3,4}:\begin{pmatrix}
e_1  &  0  &0 \\
0  &  e_2  &0\\
0  &  e_3  &0
\end{pmatrix},
\mathcal{J}_{3,5}:\begin{pmatrix}
e_1  &  0  &0 \\
0  &  e_2  &2e_3\\
0  &  -e_3  &0
\end{pmatrix},
\mathcal{J}_{3,6}:\begin{pmatrix}
e_1  &  e_2  &e_3 \\
-e_2  &  e_2  &0\\
-e_3  &  e_3  &0
\end{pmatrix},\\&&
\mathcal{J}_{3,7}:\begin{pmatrix}
e_1  &  e_2  &e_3 \\
-e_2  &  e_2  &e_3\\
-e_3  &  0  &0
\end{pmatrix},
\mathcal{J}_{3,8}:\begin{pmatrix}
e_1  &  0  &e_3 \\
0  &  e_2  &e_3\\
-e_3  &  0  &0
\end{pmatrix},
\mathcal{J}_{3,9}:\begin{pmatrix}
e_1  &  e_2  &0 \\
-e_2  &  e_2  &e_3\\
0  &  0  &0
\end{pmatrix},\\&&
\mathcal{J}_{3,10}:\begin{pmatrix}
e_1  &  0  &e_2+e_3 \\
0  &  e_2  &e_3\\
-e_2-e_3  &  0  &0
\end{pmatrix},
\mathcal{J}_{3,11}:\begin{pmatrix}
e_1  &  -e_1  &e_3 \\
e_1  &  e_2  &0\\
-e_3  &  e_3  &0
\end{pmatrix},\\&&
\mathcal{J}_{3,12}:\begin{pmatrix}
e_1  &  -e_1  &e_1+e_2+e_3 \\
e_1  &  e_2  &-e_1-e_2\\
-e_1-e_2-e_3  &  e_1+e_2+e_3  &0
\end{pmatrix},
\mathcal{J}_{3,13}:\begin{pmatrix}
e_1  &  0  &0 \\
0  &  e_2  &e_3\\
0  &  0  &0
\end{pmatrix}.\end{eqnarray*}
\\
\end{proposition}
\begin{proposition}{\rm (case $\mathcal{J}_4$)} The formal
characteristic matrix of the Jordan algebra
$(\mathcal{J}_4,\circ)$ is given by$$\begin{pmatrix}
2e_1  &  e_2  &2e_3 \\
e_2  &  0  &0\\
2e_3  &  0  &0
\end{pmatrix}.$$Any compatible pre-Jordan algebra structure on
$(\mathcal{J}_4,\circ)$ is isomorphic to one of the following
(mutually non-isomorphic) pre-Jordan
algebras given by their formal characteristic matrices respectively.
\\
\begin{eqnarray*}
&&\mathcal{J}_{4,1}(\alpha ):\begin{pmatrix}
e_1  &  (\alpha +1)e_2+\frac{1}{2}e_3  &-2\alpha (\alpha -1)e_2+(2-\alpha )e_3 \\
-\alpha e_2-\frac{1}{2}e_3  &  0  &0\\
2\alpha (\alpha -1)e_2+\alpha e_3  &  0  &0
\end{pmatrix},\\&&
\mathcal{J}_{4,2}:\begin{pmatrix}
e_1  &  0  &e_3 \\
e_2  &  0  &0\\
e_3  &  0  &0
\end{pmatrix},
\mathcal{J}_{4,3}:\begin{pmatrix}
e_1  &  0  &2e_3 \\
e_2  &  0  &0\\
0  &  0  &0
\end{pmatrix},
\mathcal{J}_{4,4}:\begin{pmatrix}
e_1  &  e_3  &e_3 \\
e_2-e_3  &  0  &0\\
e_3  &  0  &0
\end{pmatrix},\\&&
\mathcal{J}_{4,5}:\begin{pmatrix}
e_1  &  e_3  &2e_3 \\
e_2-e_3  &  0  &0\\
0  &  0  &0
\end{pmatrix},
\mathcal{J}_{4,6}:\begin{pmatrix}
e_1  &  e_2  &e_3 \\
0  &  0  &0\\
e_3  &  0  &0
\end{pmatrix},
\mathcal{J}_{4,7}:\begin{pmatrix}
e_1  &  2e_2  &2e_3 \\
-e_2  &  0  &0\\
0  &  0  &0
\end{pmatrix},\\&&
\mathcal{J}_{4,8}:\begin{pmatrix}
e_1  &  e_2  &2e_3 \\
0  &  0  &0\\
0  &  0  &0
\end{pmatrix},
\mathcal{J}_{4,9}:\begin{pmatrix}
e_1  &  e_2  &2e_3+e_2 \\
0  &  0  &0\\
-e_2  &  0  &0
\end{pmatrix},
\mathcal{J}_{4,10}:\begin{pmatrix}
e_1  &  2e_2  &e_3 \\
-e_2  &  0  &0\\
e_3  &  0  &0
\end{pmatrix},\\&&
\mathcal{J}_{4,11}:\begin{pmatrix}
e_1  &  2e_2  &e_2+e_3 \\
-e_2  &  0  &0\\
e_3-e_2  &  0  &0
\end{pmatrix},
\mathcal{J}_{4,12}:\begin{pmatrix}
e_1  &  e_2  &2e_3 \\
0  &  0  &e_1\\
0  &  -e_1  &0
\end{pmatrix},\\&&
\mathcal{J}_{4,13}:\begin{pmatrix}
e_1  &  e_2  &e_1-e_2+2e_3 \\
0  &  0  &e_1\\
-e_1+e_2  &  -e_1  &0
\end{pmatrix}.\end{eqnarray*}
\\
\end{proposition}
\begin{proposition}{\rm (case $\mathcal{J}_5$)} The formal
characteristic matrix of the Jordan algebra
$(\mathcal{J}_5,\circ)$ is given by$$\begin{pmatrix}
2e_1  &  e_2  &e_3 \\
e_2  &  0  &0\\
e_3  &  0  &0
\end{pmatrix}.$$Any compatible pre-Jordan algebra structure on
$(\mathcal{J}_5,\circ)$ is isomorphic to one of the following
(mutually non-isomorphic) pre-Jordan
algebras given by their formal characteristic matrices respectively.\\
\begin{eqnarray*}
&&
\mathcal{J}_{5,1}:\begin{pmatrix}
e_1  &  0  &0 \\
e_2  &  0  &0\\
e_3  &  0  &0
\end{pmatrix},
\mathcal{J}_{5,2}:\begin{pmatrix}
e_1  &  2e_2  &2e_3 \\
-e_2  &  0  &0\\
-e_3  &  0  &0
\end{pmatrix},
\mathcal{J}_{5,3}:\begin{pmatrix}
e_1  &  e_2  &e_3 \\
0  &  0  &0\\
0  &  0  &0
\end{pmatrix},\\&&
\mathcal{J}_{5,4}:\begin{pmatrix}
e_1  &  e_2  &0 \\
0  &  0  &0\\
e_3  &  0  &0
\end{pmatrix},
\mathcal{J}_{5,5}:\begin{pmatrix}
e_1  &  e_2  &2e_3 \\
0  &  0  &0\\
-e_3  &  0  &0
\end{pmatrix},
\mathcal{J}_{5,6}:\begin{pmatrix}
e_1  &  2e_2  &0 \\
-e_2  &  0  &0\\
e_3  &  0  &0
\end{pmatrix}.
\end{eqnarray*}
\end{proposition}
\begin{proposition}{\rm (case $\mathcal{J}_6$)} The formal
characteristic matrix of the Jordan algebra
$(\mathcal{J}_6,\circ)$ is given by$$\begin{pmatrix}
2e_1  &  e_2  &0 \\
e_2  &  2e_3  &0\\
0  &  0  &0
\end{pmatrix}.$$Any compatible pre-Jordan algebra structure on
$(\mathcal{J}_6,\circ)$ is isomorphic to one of the following
(mutually non-isomorphic) pre-Jordan
algebras given by their formal characteristic matrices respectively.
\\
\begin{eqnarray*}
\mathcal{J}_{6,1}:\begin{pmatrix}
e_1  &  2e_2  &e_3 \\
-e_2  &  e_3  &0\\
-e_3  &  0  &0
\end{pmatrix},
\mathcal{J}_{6,2}:\begin{pmatrix}
e_1  &  0  &e_3 \\
e_2  &  e_3  &0\\
-e_3  &  0  &0
\end{pmatrix},
\mathcal{J}_{6,3}:\begin{pmatrix}
e_1  &  e_2  &0 \\
0  &  e_3  &0\\
0  &  0  &0
\end{pmatrix},
\mathcal{J}_{6,4}:\begin{pmatrix}
e_1  &  e_2  &0 \\
0  &  e_3  &e_1\\
0  &  -e_1  &0
\end{pmatrix}.\end{eqnarray*}
\end{proposition}
\begin{proposition}{\rm (case $\mathcal{J}_7$)} The formal
characteristic matrix of the Jordan algebra
$(\mathcal{J}_7,\circ)$ is given by$$\begin{pmatrix}
2e_1  &  e_2  &2e_3 \\
e_2  &  2e_3  &0\\
2e_3  &  0  &0
\end{pmatrix}.$$Any compatible pre-Jordan algebra structure on
$(\mathcal{J}_7,\circ)$ is isomorphic to one of the following
(mutually non-isomorphic) pre-Jordan
algebras given by their formal characteristic matrices respectively.
\\
\begin{eqnarray*}
\mathcal{J}_{7,1}:\begin{pmatrix}
e_1  &  0  &e_3 \\
e_2  &  e_3  &0\\
e_3  &  0  &0
\end{pmatrix},
\mathcal{J}_{7,2}:\begin{pmatrix}
e_1  &  e_2  &2e_3 \\
0  &  e_3  &e_1\\
0  &  -e_1  &0
\end{pmatrix},
\mathcal{J}_{7,3}:\begin{pmatrix}
e_1  &  e_2  &2e_3 \\
0  &  e_3  &0\\
0  &  0  &0
\end{pmatrix},
\mathcal{J}_{7,4}:\begin{pmatrix}
e_1  &  2e_2  &e_3 \\
-e_2  &  e_3  &0\\
e_3  &  0  &0
\end{pmatrix}.\end{eqnarray*}
\end{proposition}
\begin{proposition}{\rm (case $\mathcal{J}_8$)} The formal
characteristic matrix of the Jordan algebra
$(\mathcal{J}_8,\circ)$ is given by$$\begin{pmatrix}
0  &  0  &0 \\
0  &  2e_2  &e_3\\
0  &  e_3  &0
\end{pmatrix}.$$Any compatible pre-Jordan algebra structure on
$(\mathcal{J}_8,\circ)$ is isomorphic to one of the following
(mutually non-isomorphic) pre-Jordan
algebras given by their formal characteristic matrices respectively.
\\
\begin{eqnarray*}
&&
\mathcal{J}_{8,1}(\alpha ):\begin{pmatrix}
0  &  \alpha e_1+e_3  &0 \\
-\alpha e_1-e_3  &  e_2  &(\alpha +1)(\alpha e_1+e_3)\\
0  &  -\alpha ((\alpha +1)e_1+e_3)  &0
\end{pmatrix},
\mathcal{J}_{8,2}:\begin{pmatrix}
0  &  -e_1  &0 \\
e_1  &  e_2  &2e_3\\
0  &  -e_3  &0
\end{pmatrix},\\&&
\mathcal{J}_{8,3}:\begin{pmatrix}
0  &  0  &0 \\
0  &  e_2  &2e_3\\
0  &  -e_3  &0
\end{pmatrix},
\mathcal{J}_{8,4}:\begin{pmatrix}
0  &  -e_1  &0 \\
e_1  &  e_2  &e_1+2e_3\\
0  &  -e_1-e_3  &0
\end{pmatrix},
\mathcal{J}_{8,5}:\begin{pmatrix}
0  &  0  &0 \\
0  &  e_2  &e_1+2e_3\\
0  &  -e_1-e_3  &0
\end{pmatrix},\\&&
\mathcal{J}_{8,6}:\begin{pmatrix}
0  &  -e_1  &0 \\
e_1  &  e_2  &e_3\\
0  &  0  &0
\end{pmatrix},
\mathcal{J}_{8,7}:\begin{pmatrix}
0  &  0  &0 \\
0  &  e_2  &0\\
0  &  e_3  &0
\end{pmatrix},
\mathcal{J}_{8,8}:\begin{pmatrix}
0  &  -e_1  &0 \\
e_1  &  e_2  &0\\
0  &  e_3  &0
\end{pmatrix},
\mathcal{J}_{8,9}:\begin{pmatrix}
0  &  0  &0 \\
0  &  e_2  &e_3\\
0  &  0  &0
\end{pmatrix},\\&&
\mathcal{J}_{8,10}:\begin{pmatrix}
0  &  0  &0 \\
0  &  e_2  &e_1+e_3\\
0  &  -e_1  &0
\end{pmatrix},
\mathcal{J}_{8,11}:\begin{pmatrix}
0  &  -e_1  &0 \\
e_1  &  e_2  &e_1\\
0  &  -e_1+e_3  &0
\end{pmatrix},
\mathcal{J}_{8,12}:\begin{pmatrix}
0  &  0  &e_2 \\
0  &  e_2  &e_3\\
-e_2  &  0  &0
\end{pmatrix}.\end{eqnarray*}
\end{proposition}
\begin{proposition}{\rm (case $A_{{\bf 1}_1}$)} The formal
characteristic matrix of the Jordan algebra
$(A_{{\bf 1}_1},\circ)$ is given by$$\begin{pmatrix}
2e_1  &  0  &0 \\
0  &  2e_2  &0\\
0  &  0  &2e_3
\end{pmatrix}.$$Any compatible pre-Jordan algebra structure on
$(A_{{\bf 1}_1},\circ)$ is isomorphic to one of the following
(mutually non-isomorphic) pre-Jordan
algebras given by their formal characteristic matrices respectively.
\begin{eqnarray*}
&&A_{{\bf 1}_1,1}:\begin{pmatrix}
e_1  &  0  &0 \\
0  &  e_2  &0\\
0  &  0  &e_3
\end{pmatrix},
A_{{\bf 1}_1,2}:\begin{pmatrix}
e_1  &  -e_1  &0 \\
e_1  &  e_2  &0\\
0  &  0  &e_3
\end{pmatrix},\\&&
A_{{\bf 1}_1,3}:\begin{pmatrix}
e_1  &  e_2  &e_3 \\
-e_2  &  e_2  &0\\
-e_3  &  0  &e_3
\end{pmatrix},
A_{{\bf 1}_1,4}:\begin{pmatrix}
e_1  &  e_2+e_3  &0 \\
-e_2-e_3  &  e_2  &e_3\\
0  &  -e_3  &e_3
\end{pmatrix}.
\end{eqnarray*}
\end{proposition}
\begin{proposition}{\rm (case $A_{{\bf 1}_2}$)} The formal
characteristic matrix of the Jordan algebra
$(A_{{\bf 1}_2},\circ)$ is given by$$\begin{pmatrix}
2e_1  &  0  &0 \\
0  &  2e_2  &2e_3\\
0  &  2e_3  &0
\end{pmatrix}.$$Any compatible pre-Jordan algebra structure on
$(A_{{\bf 1}_2},\circ)$ is isomorphic to one of the following
(mutually non-isomorphic) pre-Jordan
algebras given by their formal characteristic matrices respectively.
\begin{eqnarray*}
&&
A_{{\bf 1}_2,1}:\begin{pmatrix}
e_1  &  -e_1  &0 \\
e_1  &  e_2  &e_3\\
0  &  e_3  &0
\end{pmatrix},
A_{{\bf 1}_2,2}:\begin{pmatrix}
e_1  &  -e_1  &0 \\
e_1  &  e_2  &2e_3\\
0  &  0  &0
\end{pmatrix},
A_{{\bf 1}_2,3}:\begin{pmatrix}
e_1  &  -e_1  &e_3 \\
e_1  &  e_2  &e_3\\
-e_3  &  e_3  &0
\end{pmatrix},\\&&
A_{{\bf 1}_2,4}:\begin{pmatrix}
e_1  &  0  &0 \\
0  &  e_2  &e_3\\
0  &  e_3  &0
\end{pmatrix},
A_{{\bf 1}_2,5}:\begin{pmatrix}
e_1  &  0  &0 \\
0  &  e_2  &2e_3\\
0  &  0  &0
\end{pmatrix},
A_{{\bf 1}_2,6}:\begin{pmatrix}
e_1  &  e_2  &0 \\
-e_2  &  e_2  &2e_3\\
0  &  0  &0
\end{pmatrix},\\&&
A_{{\bf 1}_2,7}:\begin{pmatrix}
e_1  &  e_2  &e_3 \\
-e_2  &  e_2  &e_3\\
-e_3  &  e_3  &0
\end{pmatrix}.
\end{eqnarray*}
\end{proposition}

\begin{proposition}{\rm (case $A_{{\bf 1}_3}$)} The formal
characteristic matrix of the Jordan algebra
$(A_{{\bf 1}_3},\circ)$ is given by$$\begin{pmatrix}
2e_1  &  2e_2  &2e_3 \\
2e_2  &  2e_3  &0\\
2e_3  &  0  &0
\end{pmatrix}.$$Any compatible pre-Jordan algebra structure on
$(A_{{\bf 1}_3},\circ)$ is isomorphic to one of the following
(mutually non-isomorphic) pre-Jordan
algebras given by their formal characteristic matrices respectively.
\\
\begin{equation*}
A_{{\bf 1}_3,1}:\begin{pmatrix}
e_1  &  e_2  &e_3 \\
e_2  &  e_3  &0\\
e_3  &  0  &0
\end{pmatrix},
A_{{\bf 1}_3,2}:\begin{pmatrix}
e_1  &  2e_2  &2e_3 \\
0  &  e_3  &0\\
0  &  0  &0
\end{pmatrix}.\end{equation*}
\end{proposition}
\begin{proposition}{\rm (case $A_{{\bf 1}_4}$)} The formal
characteristic matrix of the Jordan algebra
$(A_{{\bf 1}_4},\circ)$ is given by$$\begin{pmatrix}
2e_1  &  2e_2  &2e_3 \\
2e_2  &  0  &0\\
2e_3  &  0  &0
\end{pmatrix}.$$Any compatible pre-Jordan algebra structure on
$(A_{{\bf 1}_4},\circ)$ is isomorphic to one of the following
(mutually non-isomorphic) pre-Jordan
algebras given by their formal characteristic matrices respectively.
\begin{equation*}
A_{{\bf 1}_4,1}:\begin{pmatrix}
e_1  &  e_2  &2e_3 \\
e_2  &  0  &0\\
0  &  0  &0
\end{pmatrix},
A_{{\bf 1}_4,2}:\begin{pmatrix}
e_1  &  e_2  &e_3 \\
e_2  &  0  &0\\
e_3  &  0  &0
\end{pmatrix},
A_{{\bf 1}_4,3}:\begin{pmatrix}
e_1  &  2e_2  &2e_3 \\
0  &  0  &0\\
0  &  0  &0
\end{pmatrix}.
\end{equation*}
\end{proposition}
\begin{proposition}{\rm (case $A_1$)} The formal
characteristic matrix of the Jordan algebra
$(A_1,\circ)$ is given by$$\begin{pmatrix}
2e_1  &  0  &0 \\
0  &  2e_2  &0\\
0  &  0  &0
\end{pmatrix}.$$Any compatible pre-Jordan algebra structure on
$(A_1,\circ)$ is isomorphic to one of the following
(mutually non-isomorphic) pre-Jordan
algebras given by their formal characteristic matrices respectively.
\begin{eqnarray*}
A_{1,1}:\begin{pmatrix}
e_1  &  0  &0 \\
0  &  e_2  &0\\
0  &  0  &0
\end{pmatrix},
A_{1,2}:\begin{pmatrix}
e_1  &  0  &e_3 \\
0  &  e_2  &0\\
e_3  &  0  &0
\end{pmatrix},
A_{1,3}:\begin{pmatrix}
e_1  &  e_2  &0 \\
-e_2  &  e_2  &0\\
0  &  0  &0
\end{pmatrix},
A_{1,4}:\begin{pmatrix}
e_1  &  e_2  &e_3 \\
-e_2  &  e_2  &0\\
-e_3  &  0  &0
\end{pmatrix}.
\end{eqnarray*}
\end{proposition}
\begin{proposition}{\rm (case $A_2$)} The formal
characteristic matrix of the Jordan algebra
$(A_2,\circ)$ is given by$$\begin{pmatrix}
0  &  0  &0 \\
0  &  2e_2  &2e_3\\
0  &  2e_3  &0
\end{pmatrix}.$$Any compatible pre-Jordan algebra structure on
$(A_2,\circ)$ is isomorphic to one of the following
(mutually non-isomorphic) pre-Jordan
algebras given by their formal characteristic matrices respectively.
\begin{eqnarray*}
&&A_{2,1}:\begin{pmatrix}
0  &  -e_1  &0 \\
e_1  &  e_2  &e_3\\
0  &  e_3  &0
\end{pmatrix},
A_{2,2}:\begin{pmatrix}
0  &  0  &0 \\
0  &  e_2  &2e_3\\
0  &  0  &0
\end{pmatrix},
A_{2,3}:\begin{pmatrix}
0  &  -e_1  &0 \\
e_1  &  e_2  &2e_3\\
0  &  0  &0
\end{pmatrix},\\&&
A_{2,4}:\begin{pmatrix}
0  &  -e_1  &0 \\
e_1  &  e_2  &e_1+2e_3\\
0  &  -e_1  &0
\end{pmatrix},
A_{2,5}:\begin{pmatrix}
0  &  0  &0 \\
0  &  e_2  &e_3\\
0  &  e_3  &0
\end{pmatrix},
A_{2,6}:\begin{pmatrix}
0  &  e_3  &0 \\
-e_3  &  e_2  &e_3\\
0  &  e_3  &0
\end{pmatrix}.
\end{eqnarray*}
\end{proposition}
\begin{proposition}{\rm (case $A_3$)} The formal
characteristic matrix of the Jordan algebra
$(A_3,\circ)$ is given by$$\begin{pmatrix}
2e_1  &  0  &0 \\
0  &  2e_3  &0\\
0  &  0  &0
\end{pmatrix}.$$Any compatible pre-Jordan algebra structure on
$(A_3,\circ)$ is isomorphic to one of the following
(mutually non-isomorphic) pre-Jordan
algebras given by their formal characteristic matrices respectively.
\begin{equation*}
A_{3,1}:\begin{pmatrix}
e_1  &  0  &0 \\
0  &  e_3  &0\\
0  &  0  &0
\end{pmatrix},
A_{3,2}:\begin{pmatrix}
e_1  &  e_2  &e_3 \\
-e_2  &  e_3  &0\\
-e_3  &  0  &0
\end{pmatrix}.
\end{equation*}
\end{proposition}
\begin{proposition}{\rm (case $A_4$)} The formal
characteristic matrix of the Jordan algebra
$(A_4,\circ)$ is given by$$\begin{pmatrix}
2e_2  &  2e_3  &0 \\
2e_3  &  0  &0\\
0  &  0  &0
\end{pmatrix}.$$Any compatible pre-Jordan algebra structure on
$(A_4,\circ)$ is isomorphic to one of the following
(mutually non-isomorphic) pre-Jordan
algebras given by their formal characteristic matrices respectively.
$$
A_{4,1}(\alpha ):\begin{pmatrix}
e_2  &  (1+\alpha)e_3  &0 \\
(1-\alpha )e_3  &  0  &0\\
0  &  0  &0
\end{pmatrix}.
$$
\end{proposition}
\begin{proposition}{\rm (case $A_5$)} The formal
characteristic matrix of the Jordan algebra
$(A_5,\circ)$ is given by$$\begin{pmatrix}
2e_1  &  0  &0 \\
0  &  0  &0\\
0  &  0  &0
\end{pmatrix}.$$Any compatible pre-Jordan algebra structure on
$(A_5,\circ)$ is isomorphic to one of the following
(mutually non-isomorphic) pre-Jordan
algebras given by their formal characteristic matrices respectively.
\begin{equation*}
A_{5,1}:\begin{pmatrix}
e_1  &  0  &0 \\
0  &  0  &0\\
0  &  0  &0
\end{pmatrix},
A_{5,2}:\begin{pmatrix}
e_1  &  e_2  &e_3 \\
-e_2  &  0  &0\\
-e_3  &  0  &0
\end{pmatrix},
A_{5,3}:\begin{pmatrix}
e_1  &  e_2  &0 \\
-e_2  &  0  &0\\
0  &  0  &0
\end{pmatrix}.
\end{equation*}
\end{proposition}
\begin{proposition}{\rm (case $A_6$)} The formal
characteristic matrix of the Jordan algebra
$(A_6,\circ)$ is given by$$\begin{pmatrix}
0  &  2e_3  &0 \\
2e_3  &  0  &0\\
0  &  0  &0
\end{pmatrix}.$$Any compatible pre-Jordan algebra structure on
$(A_6,\circ)$ is isomorphic to one of the following
(mutually non-isomorphic) pre-Jordan
algebras given by their formal characteristic matrices respectively.
\begin{equation*}
A_{6,1}(\alpha ):\begin{pmatrix}
0  &  (1+\sqrt{\alpha })e_3  &0 \\
(1-\sqrt{\alpha })e_3  &  0  &0\\
0  &  0  &0
\end{pmatrix},
A_{6,2}:\begin{pmatrix}
0  &  2e_3  &0 \\
0  &  0  &e_1\\
0  &  -e_1  &0
\end{pmatrix}.
\end{equation*}
\end{proposition}
\begin{proposition}{\rm (case $A_7$)} The formal
characteristic matrix of the Jordan algebra
$(A_7,\circ)$ is given by$$\begin{pmatrix}
0  &  0  &0 \\
0  &  2e_3  &0\\
0  &  0  &0
\end{pmatrix}.$$Any compatible pre-Jordan algebra structure on
$(A_7,\circ)$ is isomorphic to one of the following
(mutually non-isomorphic) pre-Jordan
algebras given by their formal characteristic matrices respectively.
\begin{equation*}
A_{7,1}:\begin{pmatrix}
0  &  0  &0 \\
0  &  e_3  &e_1\\
0  &  -e_1  &0
\end{pmatrix},
A_{7,2}:\begin{pmatrix}
0  &  0  &0 \\
0  &  e_3  &0\\
0  &  0  &0
\end{pmatrix},
A_{7,3}:\begin{pmatrix}
0  &  e_3  &0 \\
-e_3  &  e_3  &0\\
0  &  0  &0
\end{pmatrix}.
\end{equation*}
\end{proposition}
\begin{proposition}{\rm (case $A_8$)} The formal
characteristic matrix of the Jordan algebra
$(A_8,\circ)$ is given by$$\begin{pmatrix}
0  &  0  &0 \\
0  &  0  &0\\
0  &  0  &0
\end{pmatrix}.$$Any compatible pre-Jordan algebra structure on
$(A_8,\circ)$ is isomorphic to one of the following
(mutually non-isomorphic) pre-Jordan
algebras given by their formal characteristic matrices respectively.
\begin{equation*}
A_{8,1}:\begin{pmatrix}
0  &  0  &0 \\
0  &  0  &0\\
0  &  0  &0
\end{pmatrix},
A_{8,2}:\begin{pmatrix}
0  &  0  &0 \\
0  &  0  &e_1\\
0  &  -e_1  &0
\end{pmatrix}.
\end{equation*}
\end{proposition}
\begin{theorem}
	Every $3$-dimensional pre-Jordan algebra is isomorphic to one of the pre-Jordan algebras listed in this section.
\end{theorem}
\begin{proof}
	It follows by Lemma~\ref{iso}.
\end{proof}
\begin{corollary}
Adopting the notations of $3$-dimensional complex associative algebras in \cite{KSTT}, we have the following correspondences for 3-dimensional associative pre-Jordan algebras:\\
$\mathcal{J}_{2,4}\simeq U_1^3,\;\mathcal{J}_{3,4}\simeq W_8^3,\;\mathcal{J}_{3,13}\simeq W_7^3,\;\mathcal{J}_{4,2}\simeq W_9^3,\;\mathcal{J}_{4,6}\simeq W_{10}^3,\;\mathcal{J}_{5,1}\simeq C_4^3,\;\mathcal{J}_{5,3}\simeq C_3^3,\; \mathcal{J}_{5,4}\simeq C_2^3,\;\mathcal{J}_{8,7}\simeq W_6^3,\\ \mathcal{J}_{8,9}\simeq W_5^3,\; A_{{\bf1}_1,1}\simeq U_2^3,\; A_{{\bf1}_2,4}\simeq U_3^3,\; A_{{\bf1}_3,1}\simeq U_4^3,\; A_{{\bf1}_4,2}\simeq U_0^3,\; A_{1,1}\simeq S^3_3,\; A_{2,5}\simeq S_4^3,\; A_{3,1}\simeq S_2^3,\\ A_{4,1}(\alpha=0)\simeq S_1^3,\; A_{5,1}\simeq W_4^3,\; A_{6,1}(\alpha=1)\simeq W_2^3,\; A_{6,1}(\alpha=\frac{k^2}{k^2-4})\simeq W^3_3(k),\; A_{7,2}\simeq W_1^3,\; A_{7,3}\simeq W_3^3(2),\\ A_{8,1}\simeq C_0^3,\; A_{8,2}\simeq C_1^3$.
\end{corollary}

\section{Rota-Baxter operators on Jordan algebras in dimensions $\leq 3$
and the induced pre-Jordan algebras}

In this section, we give all Rota-Baxter operators on complex Jordan algebras in dimensions $\leq 3$ and the induced pre-Jordan algebras. If not assigned concretely, the parameters can be any complex numbers.

\begin{definition}\label{Rota-Baxter} {\rm (\cite{HNB})
    A {\em Rota-Baxter operator (of weight 0)} on a Jordan algebra $(J,\circ)$ is a linear transformation $R:J\rightarrow J$ satisfying the following identity:
\begin{equation}
R(x)\circ R(y) = R(R(x)\circ y + x\circ R(y)),\;\;\forall x,y\in
J.
\end{equation}}
\end{definition}

\begin{proposition} {\rm (\cite{HNB})} Let $(J,\circ)$ be a Jordan
algebra and $R:J\rightarrow J$ be a Rota-Baxter operator. Then the
bilinear multiplication $\cdot:J\times J\rightarrow J$ given by
\begin{equation}\label{eq:induceJ}
    x\cdot y := R(x)\circ y,\;\;\forall x,y\in J,
    \end{equation}
defines a pre-Jordan algebra.
\end{proposition}

\begin{proposition}
Keep the notations as in Proposition~\ref{prop:1dim}. The
Rota-Baxter operators on the Jordan algebra $(A,\circ_1)$ given by
$e\circ e=e$ are zero and  the Rota-Baxter operators on the Jordan
algebra $(A,\circ_2)$ given by $e\circ e=0$ are all elements of ${\rm End}(\complex)$. Both of them induce the pre-Jordan algebra $(\complex,\cdot)$
given by $e\cdot e=0$.
\end{proposition}

\begin{proof}
Note the first conclusion can be obtained from
the fact that all Rota-Baxter operators on a unital finite
dimensional complex associative algebra are nilpotent (\cite{G}).
\end{proof}

Let $(J,\circ)$ be a Jordan algebra with a basis
$\{e_1,\cdots,e_N\}$. Let $R:J\rightarrow J$ be a linear
transformation on $J$. Set
\begin{equation}S(x,y)=R(x)\circ R(y)- R(R(x)\circ y + R(y)\circ
x),\;\;\forall x,y\in J,\end{equation}
\begin{equation}
e_i\circ e_j=\sum_{k=1}^n c_{ij}^k e_k,\;\;R(e_i)=\sum_{j=1}^N
r_{ji}e_j.
\end{equation}
Hence we have \begin{equation}S(e_j,e_k)=\sum_{i=1}^N
S^i_{jk}e_i,\end{equation} where
\begin{equation}
S_{jk}^i=\sum_{n=1}^{N}\sum_{m=1}^{N}\left(r_{nj}r_{mk}c_{mn}^i-r_{nj}r_{im}c_{nk}^m-r_{nk}r_{im}c_{nj}^m\right).
\end{equation}
The following conclusion follows immediately.

\begin{proposition}
Keep the notations as above and let $R=\left(\begin{matrix} r_{11}
&\cdots &r_{1n}\cr \cdots & \cdots &\cdots\cr r_{n1}& \cdots
&r_{nn}\cr \end{matrix}\right)$. Then it is a Rota-Baxter operator on the
Jordan algebra $(J,\circ)$ if and only if the entries $r_{ij}$ satisfy the set of
equations $\{S^k_{ij}=0\}$.
\end{proposition}
\begin{proposition}
    Any Rota-Baxter operator on the Jordan algebra $(J_1,\circ)$ is zero, which induces the pre-Jordan algebra isomorphic to $J_{6,1}$ through Eq.~(\ref{eq:induceJ}).
\end{proposition}
\begin{proof}
There are still the following ``essentially available" equations
in the set of equations $\{S^k_{ij}=0\}$:
\begin{eqnarray*}
S^1_{11}=-2r^2_{11}=0,\;\; S^1_{22}=2r_{12}^2-4r_{12}r_{22}=0,\;\;
S^2_{11}=2r_{21}^2-4r_{11}r_{21}=0,\;\; S^2_{22}=-2r^2_{22}=0.
\end{eqnarray*}
There is only one solution:
$$r_{11}=r_{12}=r_{21}=r_{22}=0.$$
Thus the conclusion holds.
\end{proof}
\begin{proposition}\label{R2}
 Any Rota-Baxter operator on the Jordan algebra $(J_2,\circ)$ is of the form $\begin{pmatrix}

                                                       0 & 0 \\
                                                       r_{21} & 0
                                                     \end{pmatrix}$.
The induced pre-Jordan algebra
through Eq.~(\ref{eq:induceJ}) is isomorphic to $J_{6,1}$
when $r_{21}=0$ or $J_{5,1}$ otherwise, respectively.
\end{proposition}
\begin{proof}
The ``essentially available" equations in the set of equations
$\{S^k_{ij}=0\}$ are:
$$S^1_{11}=-2r_{11}^2-4r_{12}r_{21}=0,\;\;S^1_{22}=-2r_{12}^2=0,\;\;S^2_{21}=-2r^2_{22}=0.$$
The solutions are
$$r_{11}=r_{12}=r_{22}=0,\;\;r_{21}\in \complex.$$
The induced pre-Jordan algebras through Eq.~(\ref{eq:induceJ}) are
given by their formal characteristic matrices
$\left(\begin{matrix} r_{21}e_2 & 0\cr 0&0\cr\end{matrix}\right)$.
When $r_{21}=0$, it is isomorphic to the pre-Jordan algebra $J_{6,1}$.
When $r_{21}\ne 0$, it is isomorphic to $J_{5,1}$.
\end{proof}
\begin{proposition}
    Any Rota-Baxter operator on the Jordan algebra $(J_3,\circ)$ is of the form $\begin{pmatrix}
    0 & 0 \\
    r_{21} & r_{22}
    \end{pmatrix}$. The induced pre-Jordan algebra is isomorphic to $J_{6,1}$.
\end{proposition}
\begin{proof}
The ``essentially available" equations in the set of equations
$\{S^k_{ij}=0\}$ are:
    $$S^1_{11}=-2r_{11}^2=0,\;\;S^1_{22}=2r_{12}^2=0.$$
    The solutions are $$
    r_{11}=r_{12}=0,\;\;r_{21}\in\complex,\;\;r_{22}\in\complex.
    $$
    All induced multiplications vanish.
\end{proof}
\begin{proposition}
    Any Rota-Baxter operator on the Jordan algebra $(J_4,\circ)$ is of the form $\begin{pmatrix}
    \lambda\mu & -\lambda^2 \\
    \mu^2 & -\lambda\mu
    \end{pmatrix}$,\\
when $\lambda\neq0$, the induced pre-Jordan algebra is isomorphic to $J_{2,2}$,\\
when $\lambda=0,\mu\neq0$, the induced pre-Jordan algebra is isomorphic to $J_{5,1}$,\\
when $\lambda=\mu=0$, the induced pre-Jordan algebra is isomorphic to $J_{6,1}$.
\end{proposition}
\begin{proof}
    Set $r_{12}=-\lambda^2,r_{21}=\mu^2$, and then the ``essentially available" equations in the set of equations
$\{S^k_{ij}=0\}$ are:
$$
    S^1_{11}=-2r_{11}^2+2\lambda^2\mu^2=0,\;\;
    S^2_{12}=-r_{22}^2+\lambda^2\mu^2=0,\;\;
    S^1_{12}=\lambda^2(r_{11}+r_{22})=0.
$$
    By choosing the sign of $\lambda$ properly, we may assume $r_{11}=\lambda\mu$. If $r_{22}=\lambda\mu$, then $\lambda^3\mu=0,r_{22}=0=-\lambda\mu$, so $r_{22}=-\lambda\mu$ always holds.
    \newline We obtain the algebra with the following formal characteristic matrix:
$\begin{pmatrix}
  2\lambda\mu e_1+\mu^2e_2 & \lambda\mu e_2 \\
  -2\lambda^2e_1-\lambda\mu e_2 & -\lambda^2e_2
\end{pmatrix}$.
If $B\neq 0$: Under the basis $e_1^\prime = -\lambda^{-2}e_2,e_2^\prime =\lambda e_1+\mu e_2$, the algebra above has the following formal characteristic matrix:
$\begin{pmatrix}
  e_1^\prime & 2e_2^\prime \\
  0 & 0
\end{pmatrix}.$
    \newline If $\lambda = 0$: The formal characteristic matrix is $\begin{pmatrix}
                                                                    \mu^2e_2 & 0 \\
                                                                    0 & 0
                                                                  \end{pmatrix}$,
     and the rest of proof is the same as Theorem~\ref{R2}.
\end{proof}
\begin{proposition}
    Any Rota-Baxter operator on the Jordan algebra $(J_5,\circ)$ is of the following forms:\\
    \begin{enumerate}
      \item $\begin{pmatrix}
    0 & 0 \\
    r_{21} & r_{22}
    \end{pmatrix}$,\\
the induced pre-Jordan algebra is isomorphic to $J_{6,1}$;
      \item $\begin{pmatrix}
    2r_{22} & 0 \\
    r_{21} & r_{22}
    \end{pmatrix}(r_{22}\neq0)$,\\
the induced pre-Jordan algebra is isomorphic to $J_{5,1}$.
    \end{enumerate}
\end{proposition}
\begin{proof}
The ``essentially available" equations in the set of equations
$\{S^k_{ij}=0\}$ are:
$$
S^1_{12}=-2r_{12}^2,\;\;S^2_{11}=2r_{11}(r_{11}-2r_{22}).
$$
The solutions are $$
r_{11}=0\mbox{ or }2r_{22},\;\;r_{12}=0,\;\;r_{21}\in\complex,\;\;r_{22}\in\complex.
$$
    This completes the first assertion.
    \newline Pre-Jordan algebras induced by the former ones are isomorphic to $J_{6,1}$. And the latter ones induce the formal characteristic matrices $\begin{pmatrix}
                                                                                                  2r_{22}e_2 & 0 \\
                                                                                                  0 & 0
                                                                                                \end{pmatrix}$, which are isomorphic to $J_{5,1}$.
\end{proof}
\begin{proposition}
    Every linear transformation on $(J_6, \circ)$ is a Rota-Baxter operator. The induced pre-Jordan algebra is isomorphic to $J_{6,1}$.
\end{proposition}
\begin{proof}
    Obvious.
\end{proof}
\indent Summarizing the above study, we have the following conclusion:
\begin{theorem}\label{2dimcompare}
Each pre-Jordan algebra obtained from Rota-Baxter
operators on $2$-dimensional Jordan algebras through
Eq.~(\ref{eq:induceJ}) is isomorphic to one of the pre-Jordan algebras
$J_{2,2}$, $J_{5,1}$ and $J_{6,1}$.
\end{theorem}

Next we give all Rota-Baxter operators on 3-dimensional Jordan
algebras and the induced pre-Jordan algebras through
Eq.~(\ref{eq:induceJ}). We only list the results while omitting the
proof.

\begin{proposition}
Any Rota-Baxter operator on the Jordan algebra $(\mathcal{J}_1, \circ)$ is of the following forms:\\
\begin{enumerate}
  \item $\begin{pmatrix}
  0 & -\frac{r_{13}^3}{2r_{23}^2} & r_{13} \\
  0 & -\frac{r_{13}^2}{2r_{23}} & r_{23} \\
  0 & -\frac{r_{13}^4}{4r_{23}^3} & \frac{r_{13}^2}{2r_{23}}
\end{pmatrix}(r_{23}\neq0)$,\\
when $r_{13}=0$, the induced pre-Jordan algebra is isomorphic to $A_{4,1}(\alpha =1)$,\\
when $r_{13}\neq0$, the induced pre-Jordan algebra is isomorphic to $A_{2,2}$;
  \item $\begin{pmatrix}
   0 & 0 & 0 \\
   0 & 0 & 0 \\
   0 & r_{32} & 0
 \end{pmatrix}$,\\
when $r_{32}=0$, the induced pre-Jordan algebra is isomorphic to $A_{8,1}$,\\
when $r_{32}\neq0$, the induced pre-Jordan algebra is isomorphic to $A_{4,1}(\alpha =1)$.
\end{enumerate}
\end{proposition}
\begin{proposition}
Any Rota-Baxter operator on the Jordan algebra $(\mathcal{J}_2, \circ)$ is of the following forms:\\
\begin{enumerate}
  \item $\begin{pmatrix}
   0 & 0 & 0 \\
   0 & 0 & 0 \\
   r_{31} & r_{32} & 0
 \end{pmatrix}$,\\
when $r_{31}\neq r_{32}$, the induced pre-Jordan algebra is isomorphic to $A_{6,1}(\alpha =1)$,\\
when $r_{31}=r_{32}=0$, the induced pre-Jordan algebra is isomorphic to $A_{8,1}$,\\
when $r_{31}=r_{32}\neq0$, the induced pre-Jordan algebra is isomorphic to $A_{7,2}$;
  \item $\begin{pmatrix}
   0 & 0 & 0 \\
   r_{33} & -r_{33} & r_{23} \\
   \frac{r_{33}^2}{r_{23}} & \frac{-r_{33}^2}{r_{23}} & r_{33}
 \end{pmatrix}(r_{23}\neq0)$,\\
the induced pre-Jordan algebra is isomorphic to $A_{2,2}$;
  \item $\begin{pmatrix}
   -r_{33} & r_{33} & r_{13} \\
   0 & 0 & 0 \\
   -\frac{r_{33}^2}{r_{13}} & \frac{r_{33}^2}{r_{13}} & r_{33}
 \end{pmatrix}(r_{13}\neq0)$,\\
the induced pre-Jordan algebra is isomorphic to $A_{2,2}$.
\end{enumerate}
\end{proposition}
\begin{proposition}
Any Rota-Baxter operator on the Jordan algebra $(\mathcal{J}_3, \circ)$ is of the following forms:\\
\begin{enumerate}
  \item $\begin{pmatrix}
   0 & 0 & 0 \\
   0 & 0 & 0 \\
   r_{31} & r_{32} & 0
 \end{pmatrix}(r_{31}\neq0)$,\\
the induced pre-Jordan algebra is isomorphic to $A_{6,1}(\alpha =1)$;
  \item $\begin{pmatrix}
    0 & 0 & 0 \\
    0 & r_{22} & r_{23} \\
    0 & -\frac{r_{22}^2}{r_{23}} & -r_{22}
  \end{pmatrix}(r_{23}\neq0)$,\\
the induced pre-Jordan algebra is isomorphic to $A_{2,2}$;
  \item $\begin{pmatrix}
     0 & r_{12} & r_{13} \\
     0 & r_{12} & r_{13} \\
     0 & r_{32} & -r_{12}
   \end{pmatrix}(r_{12}^2+r_{13}r_{32}=0)$,\\
when $r_{12}\neq0\mbox{ or }r_{12}=r_{32}=0,r_{13}\neq0$, the induced pre-Jordan algebra is isomorphic to $A_{{\bf 1}_4,3}$,\\
when $r_{12}=0,r_{32}\neq0$, the induced pre-Jordan algebra is isomorphic to $A_{7,2}$,\\
when $r_{12}=r_{13}=r_{32}=0$, the induced pre-Jordan algebra is isomorphic to $A_{8,1}$.
\end{enumerate}
\end{proposition}
\begin{proposition}
Any Rota-Baxter operator on the Jordan algebra $(\mathcal{J}_4, \circ)$ is of the following forms:\\
\begin{enumerate}
  \item $\begin{pmatrix}
   r_{11} & r_{12} & 0 \\
   r_{21} & -r_{11} & 0 \\
   r_{31} & 0 & 0
 \end{pmatrix}(r_{11}^2+r_{12}r_{21}=0)$,\\
when $r_{31}=0$, the induced pre-Jordan algebra is isomorphic to $A_{{\bf 1}_4,3}$,\\
when $r_{31}\neq0$, the induced pre-Jordan algebra is isomorphic to $A_{{\bf 1}_3,2}$;
  \item $\begin{pmatrix}
    0 & 0 & 0 \\
    r_{21} & 0 & r_{23} \\
    0 & 0 & 0
  \end{pmatrix}(r_{23}\neq0)$,\\
the induced pre-Jordan algebra is isomorphic to $A_{6,1}(\alpha =1)$;
  \item $\begin{pmatrix}
     0 & 0 & 0 \\
     r_{21} & 0 & 0 \\
     r_{31} & 0 & 0
   \end{pmatrix}(r_{21}\neq0)$,\\
the induced pre-Jordan algebra is isomorphic to $A_{7,2}$;
  \item $\begin{pmatrix}
      0 & 0 & 0 \\
      0 & 0 & 0 \\
      r_{31} & r_{32} & 0
    \end{pmatrix}$,\\
when $r_{31}=r_{32}=0$, the induced pre-Jordan algebra is isomorphic to $A_{8,1}$,\\
when $r_{32}\neq0$, the induced pre-Jordan algebra is isomorphic to $A_{6,1}(\alpha =1)$,\\
when $r_{32}=0,r_{31}\neq0$, the induced pre-Jordan algebra is isomorphic to $A_{7,2}$;
  \item $\begin{pmatrix}
       0 & 0 & 0 \\
       \frac{r_{23}r_{31}}{r_{33}} & -2r_{33} & r_{23} \\
       r_{31} & -\frac{2r_{33}^2}{r_{23}} & r_{33}
     \end{pmatrix}(r_{23},r_{33}\neq0)$,\\
the induced pre-Jordan algebra is isomorphic to $A_{6,1}(\alpha=1)$.
\end{enumerate}
\end{proposition}
\begin{proposition}

Any Rota-Baxter operator on the Jordan algebra $(\mathcal{J}_5, \circ)$ is of the following forms:\\
\begin{enumerate}
  \item $\begin{pmatrix}
    0 & 0 & 0 \\
    0 & 0 & 0 \\
    0 & 0 & 0
  \end{pmatrix}$,\\
the induced pre-Jordan algebra is isomorphic to $A_{8,1}$;
  \item $A\begin{pmatrix}
    0 & 1 & 0 \\
    0 & 0 & 0 \\
    0 & 0 & 0
  \end{pmatrix}A^{-1}(A\in{\rm Aut}(\mathcal{J}_5,\circ))$,\\
the induced pre-Jordan algebra is isomorphic to $\mathcal{J}_{4,8}$;
  \item $A\begin{pmatrix}
    0 & 0 & 0 \\
    0 & 0 & 1 \\
    0 & 0 & 0
  \end{pmatrix}A^{-1}(A\in{\rm Aut}(\mathcal{J}_5,\circ))$,\\
the induced pre-Jordan algebra is isomorphic to $A_{6,1}(\alpha =1)$;
  \item $A\begin{pmatrix}
    0 & 0 & 0 \\
    1 & 0 & 0 \\
    0 & 0 & 0
  \end{pmatrix}A^{-1}(A\in{\rm Aut}(\mathcal{J}_5,\circ))$,\\
the induced pre-Jordan algebra is isomorphic to $A_{7,2}$, \\
where ${\rm Aut}(\mathcal{J}_5,\circ)=\begin{pmatrix}
	1 & 0 &0  \\
	\lambda_1 & \lambda_2 & \lambda_3\\
	\lambda_4 & \lambda_5 & \lambda_6
\end{pmatrix}(\lambda_2\lambda_6-\lambda_3\lambda_5\neq0)$.
\end{enumerate}
\end{proposition}
\begin{proposition}
Any Rota-Baxter operator on the Jordan algebra $(\mathcal{J}_6, \circ)$ is of the following forms:\\
\begin{enumerate}
  \item $\begin{pmatrix}
   0 & r_{12} & 0 \\
   0 & 0 & 0 \\
   0 & 0 & r_{33}
 \end{pmatrix}(r_{12}\neq0)$,\\
the induced pre-Jordan algebra is isomorphic to $A_{2,2}$;
  \item $\begin{pmatrix}
    0 & 0 & 0 \\
    0 & 0 & 0 \\
    r_{31} & r_{32} & r_{33}
  \end{pmatrix}$,\\
the induced pre-Jordan algebra is isomorphic to $A_{8,1}$;
  \item $\begin{pmatrix}
     r_{11} & -\frac{r_{11}^2}{r_{21}} & 0 \\
     r_{21} & -r_{11} & 0 \\
     r_{31} & r_{21}-\frac{2r_{11}r_{31}}{r_{21}} & -r_{11}+\frac{r_{11}^2r_{31}}{r_{21}^2}
   \end{pmatrix}(r_{21}\neq0)$,\\
when $r_{11}=0$, the induced pre-Jordan algebra is isomorphic to $A_{4,1}(\alpha =1)$,\\
when $r_{11}\neq0$, the induced pre-Jordan algebra is isomorphic to $A_{2,2}$.
\end{enumerate}
\end{proposition}
\begin{proposition}
Any Rota-Baxter operator on the Jordan algebra $(\mathcal{J}_7, \circ)$ is of the following forms:\\
\begin{enumerate}
  \item $\begin{pmatrix}
   0 & 0 & 0 \\
   0 & 0 & 0 \\
   r_{31} & r_{32} & 0
 \end{pmatrix}(r_{32}\neq0)$,\\
the induced pre-Jordan algebra is isomorphic to $A_{6,1}(\alpha =1)$;
  \item $\begin{pmatrix}
    0 & r_{12} & 0 \\
    0 & 0 & 0 \\
    r_{31} & 0 & 0
  \end{pmatrix}$,\\
when $r_{12}=r_{31}=0$, the induced pre-Jordan algebra is isomorphic to $A_{8,1}$,\\
when $r_{12}=0,r_{31}\neq0$, the induced pre-Jordan algebra is isomorphic to $A_{7,2}$,\\
when $r_{12}\neq0,r_{31}=0$, the induced pre-Jordan algebra is isomorphic to $A_{{\bf 1}_4,3}$,\\
when $r_{12}\neq0,r_{31}\neq0$, the induced pre-Jordan algebra is isomorphic to $A_{{\bf 1}_3,2}$;
  \item $\begin{pmatrix}
     r_{11} & -\frac{r_{11}^2}{r_{21}} & 0 \\
     r_{21} & -r_{11} & 0 \\
     r_{31} & r_{21} & 0
   \end{pmatrix}(r_{21}\neq0)$,\\
when $r_{11}=0$, the induced pre-Jordan algebra is isomorphic to $A_{4,1}(\alpha =0)$,\\
when $r_{11}\neq0,r_{21}^2+r_{11}r_{31}\neq0$, the induced pre-Jordan algebra is isomorphic to $A_{{\bf 1}_3,2}$,\\
when $r_{11}\neq0,r_{21}^2+r_{11}r_{31}=0$, the induced pre-Jordan algebra is isomorphic to $A_{{\bf 1}_4,3}$.
\end{enumerate}
\end{proposition}
\begin{proposition}
Any Rota-Baxter operator on the Jordan algebra $(\mathcal{J}_8, \circ)$ is of the following forms:\\
\begin{enumerate}
  \item $\begin{pmatrix}
   r_{11} & r_{12} & 0 \\
   0 & 0 & 0 \\
   r_{31} & r_{32} & 0
 \end{pmatrix}$,\\
when $r_{31}\neq0$, the induced pre-Jordan algebra is isomorphic to $A_{6,1}(\alpha =1)$,\\
when $r_{31}=0,r_{32}\neq0$, the induced pre-Jordan algebra is isomorphic to $A_{7,2}$,\\
when $r_{31}=r_{32}=0$, the induced pre-Jordan algebra is isomorphic to $A_{8,1}$;
  \item $\begin{pmatrix}
    r_{11} & 0 & 0 \\
    0 & r_{22} & r_{23} \\
    0 & -\frac{r_{22}^2}{r_{23}} & -r_{22}
  \end{pmatrix}(r_{23}\neq0)$,\\
the induced pre-Jordan algebra is isomorphic to $A_{2,2}$;
  \item $\begin{pmatrix}
     r_{11} & r_{12} & r_{13} \\
     0 & 0 & 0 \\
     0 & 0 & 0
   \end{pmatrix}(r_{13}\neq0)$,\\
the induced pre-Jordan algebra is isomorphic to $A_{8,1}$.
\end{enumerate}
\end{proposition}
\begin{proposition}
Any Rota-Baxter operator on $A_{{\bf 1}_1}$ is zero and the induced pre-Jordan algebra is isomorphic to $A_{8,1}$.
\end{proposition}
\begin{proposition}
Any Rota-Baxter operator on the Jordan algebra $(A_{{\bf 1}_2}, \circ)$ is of the following form:\\
$\begin{pmatrix}
   0 & 0 & 0 \\
   0 & 0 & 0 \\
   r_{31} & r_{32} & 0
 \end{pmatrix}$,\\
when $r_{31}\neq0$, the induced pre-Jordan algebra is isomorphic to $A_{6,1}(\alpha =1)$,\\
when $r_{31}=0,r_{32}\neq0$, the induced pre-Jordan algebra is isomorphic to $A_{7,2}$,\\
when $r_{31}=r_{32}=0$, the induced pre-Jordan algebra is isomorphic to $A_{8,1}$.
\end{proposition}
\begin{proposition}
Any Rota-Baxter operator on the Jordan algebra $(A_{{\bf 1}_3}, \circ)$ is of the following forms:\\
\begin{enumerate}
  \item $\begin{pmatrix}
   0 & 0 & 0 \\
   2r_{32} & 0 & 0 \\
   r_{31} & r_{32} & 0
 \end{pmatrix}(r_{32}\neq0)$,\\
the induced pre-Jordan algebra is isomorphic to $A_{4,1}(\alpha =\frac{1}{3})$;
  \item $\begin{pmatrix}
    0 & 0 & 0 \\
    0 & 0 & 0 \\
    r_{31} & r_{32} & 0
  \end{pmatrix}$,\\
when $r_{32}\neq0$, the induced pre-Jordan algebra is isomorphic to $A_{6,1}(\alpha =1)$,\\
when $r_{32}=0,r_{31}\neq0$, the induced pre-Jordan algebra is isomorphic to $A_{7,2}$,\\
when $r_{31}=r_{32}=0$, the induced pre-Jordan algebra is isomorphic to $A_{8,1}$.
\end{enumerate}
\end{proposition}
\begin{proposition}
Any Rota-Baxter operator on the Jordan algebra $(A_{{\bf 1}_4}, \circ)$ is of the following forms:\\
\begin{enumerate}
  \item $\begin{pmatrix}
   0 & 0 & 0 \\
   r_{21} & 0 & 0 \\
   r_{31} & 0 & 0
 \end{pmatrix}$,\\
when $r_{21}=r_{31}=0$, the induced pre-Jordan algebra is isomorphic to $A_{8,1}$,\\
otherwise, the induced pre-Jordan algebra is isomorphic to $A_{7,2}$;
  \item $A\begin{pmatrix}
    0 & 0 & 0 \\
    0 & 0 & 0 \\
    r_{31} & 1 & 0
  \end{pmatrix}A^{-1}(A\in{\rm Aut}(A_{1_4},\circ))$,\\
the induced pre-Jordan algebra is isomorphic to $A_{6,1}(\alpha =1)$,\\
where ${\rm Aut}(A_{1_4}, \circ)=\begin{pmatrix}
	1 & 0 &0\\
	0 & \lambda_1 &\lambda_2\\
	0 & \lambda_3 &\lambda_4
\end{pmatrix}(\lambda_1\lambda_4-\lambda_2\lambda_3\neq0)$.\\
\end{enumerate}
\end{proposition}
\begin{proposition}
Any Rota-Baxter operator on the Jordan algebra $(A_1, \circ)$ is of the following form:\\
$\begin{pmatrix}
   0 & 0 & 0 \\
   0 & 0 & 0 \\
   r_{31} & r_{32} & r_{33}
 \end{pmatrix}$,\\
the induced pre-Jordan algebra is isomorphic to $A_{8,1}$.
\end{proposition}
\begin{proposition}
Any Rota-Baxter operator on the Jordan algebra $(A_2, \circ)$ is of the following forms:\\
\begin{enumerate}
  \item $\begin{pmatrix}
   r_{11} & r_{12} & 0 \\
   0 & 0 & 0 \\
   r_{31} & r_{32} & 0
 \end{pmatrix}$,\\
when $r_{31}\neq0$, the induced pre-Jordan algebra is isomorphic to $A_{6,1}(\alpha =1)$,\\
when $r_{31}=0,r_{32}\neq0$, the induced pre-Jordan algebra is isomorphic to $A_{7,2}$,\\
when $r_{31}=r_{32}=0$, the induced pre-Jordan algebra is isomorphic to $A_{8,1}$;
  \item $\begin{pmatrix}
    r_{11} & r_{12} & r_{13} \\
    0 & 0 & 0 \\
    0 & 0 & 0
  \end{pmatrix}(r_{13}\neq0)$,\\
the induced pre-Jordan algebra is isomorphic to $A_{8,1}$.
\end{enumerate}
\end{proposition}
\begin{proposition}
Any Rota-Baxter operator on the Jordan algebra $(A_3, \circ)$ is of the following forms:\\
\begin{enumerate}
  \item $\begin{pmatrix}
   0 & 0 & 0 \\
   0 & 0 & 0 \\
   r_{31} & r_{32} & r_{33}
 \end{pmatrix}$,\\
 the induced pre-Jordan algebra is isomorphic to $A_{8,1}$;
  \item $\begin{pmatrix}
    0 & 0 & 0 \\
    0 & 2r_{33} & 0 \\
    r_{31} & r_{32} & r_{33}
  \end{pmatrix}(r_{33}\neq0)$,\\
  the induced pre-Jordan algebra is isomorphic to $A_{7,2}$.
\end{enumerate}
\end{proposition}
\begin{proposition}
Any Rota-Baxter operator on the Jordan algebra $(A_4, \circ)$ is of the following forms:\\
\begin{enumerate}
  \item $\begin{pmatrix}
   0 & 0 & 0 \\
   0 & 0 & 0 \\
   r_{31} & r_{32} & r_{33}
 \end{pmatrix}(r_{33}\neq0)$,\\
 the induced pre-Jordan algebra is isomorphic to $A_{8,1}$;
  \item $\begin{pmatrix}
    0 & 0 & 0 \\
    r_{21} & r_{22} & 0 \\
    r_{31} & r_{32} & 0
  \end{pmatrix}$,\\
when $r_{22}\neq0$, the induced pre-Jordan algebra is isomorphic to $A_{6,1}(\alpha =1)$,\\
when $r_{22}=0,r_{21}\neq0$, the induced pre-Jordan algebra is isomorphic to $A_{7,2}$,\\
when $r_{21}=r_{22}=0$, the induced pre-Jordan algebra is isomorphic to $A_{8,1}$;
  \item $\begin{pmatrix}
     r_{11} & 0 & 0 \\
     r_{21} & \frac{1}{2}r_{11} & 0 \\
     r_{31} & \frac{2}{3}r_{21} & \frac{1}{3}r_{11}
   \end{pmatrix}(r_{11}\neq0)$,\\
the induced pre-Jordan algebra is isomorphic to $A_{4,1}(\alpha=\frac{1}{3})$.
\end{enumerate}
\end{proposition}
\begin{proposition}
Any Rota-Baxter operator on the Jordan algebra $(A_5, \circ)$ is of the following form:\\
$\begin{pmatrix}
   0 & 0 & 0 \\
   r_{21} & r_{22} & r_{23} \\
   r_{31} & r_{32} & r_{33}
 \end{pmatrix}$,\\
the induced pre-Jordan algebra is isomorphic to $A_{8,1}$.
\end{proposition}
\begin{proposition}
Any Rota-Baxter operator on the Jordan algebra $(A_6, \circ)$ is of the following forms:\\
\begin{enumerate}
  \item $\begin{pmatrix}
   r_{11} & r_{12} & 0 \\
   0 & 0 & 0 \\
   r_{31} & r_{32} & 0
 \end{pmatrix}(r_{11}\neq0)$,\\
 the induced pre-Jordan algebra is isomorphic to $A_{6,1}(\alpha =1)$;
  \item $\begin{pmatrix}
    0 & 0 & 0 \\
    0 & 0 & 0 \\
    r_{31} & r_{32} & r_{33}
  \end{pmatrix}(r_{33}\neq0)$,\\
  the induced pre-Jordan algebra is isomorphic to $A_{8,1}$;
  \item $\begin{pmatrix}
     r_{11} & 0 & 0 \\
     0 & \frac{r_{11}r_{33}}{r_{11}-r_{33}} & 0 \\
     r_{31} & r_{32} & r_{33}
   \end{pmatrix}(r_{11}\neq0,r_{33}\neq0,r_{11}\neq r_{33})$,\\
   the induced pre-Jordan algebra is isomorphic to $A_{6,1}(\alpha =(1-\frac{2r_{33}}{r_{11}})^2)$;
  \item $\begin{pmatrix}
      0 & r_{12} & 0 \\
      0 & 0 & 0 \\
      r_{31} & r_{32} & 0
    \end{pmatrix}$,\\
    when $r_{12}\neq0$, the induced pre-Jordan algebra is isomorphic to $A_{7,2}$,\\
    when $r_{12}=0$, the induced pre-Jordan algebra is isomorphic to $A_{8,1}$;
  \item $\begin{pmatrix}
       r_{11} & \frac{r_{11}^2}{r_{21}} & 0 \\
       r_{21} & r_{11} & 0 \\
       r_{31} & r_{32} & r_{11}
     \end{pmatrix}(r_{21}\neq0)$,\\
     the induced pre-Jordan algebra is isomorphic to $A_{7,2}$;
  \item $\begin{pmatrix}
        0 & 0 & 0 \\
        r_{21} & r_{22} & 0 \\
        r_{31} & r_{32} & 0
      \end{pmatrix}(r_{22}\neq0)$,\\
      the induced pre-Jordan algebra is isomorphic to $A_{6,1}(\alpha =1)$.
\end{enumerate}
\end{proposition}
\begin{proposition}
Any Rota-Baxter operator on the Jordan algebra $(A_7, \circ)$ is of the following forms:\\
\begin{enumerate}
  \item $\begin{pmatrix}
   r_{11} & r_{12} & r_{13} \\
   0 & 0 & 0 \\
   r_{31} & r_{32} & r_{33}
 \end{pmatrix}$,\\
 the induced pre-Jordan algebra is isomorphic to $A_{8,1}$;
  \item $\begin{pmatrix}
    r_{11} & r_{12} & 0 \\
    0 & 2r_{33} & 0 \\
    r_{31} & r_{32} & r_{33}
  \end{pmatrix}(r_{33}\neq0)$,\\
  the induced pre-Jordan algebra is isomorphic to $A_{7,2}$.
\end{enumerate}
\end{proposition}
\begin{proposition}
All linear transformations of $A_8$ are Rota-Baxter operators. The induced pre-Jordan algebra is isomorphic to $A_{8,1}$.
\end{proposition}
Summarizing the above study, we have the following conclusion:
\begin{theorem}\label{3dimcompare}
Each $3$-dimensional pre-Jordan algebra induced by a Rota-Baxter operator on a Jordan algebra is isomorphic to one of the following pre-Jordan algebras:\\ $\mathcal{J}_{4,8},\;A_{{\bf 1}_3,2},\;A_{{\bf 1}_4,3},\;A_{2,2},\;A_{4,1}(\alpha =0,\frac{1}{3},1),\;A_{6,1}(\alpha \in\complex),\;A_{7,2},\;A_{8,1}$.
\end{theorem}
\smallskip

 \noindent
 {\bf Acknowledgements.}  The authors thank Chengming Bai for
 guidance and important discussion. All authors are supported by Innovation and Entrepreneurship Training Program for College Students of Tianjin 202010055305.

\end{document}